\documentclass[11pt]{article}

\usepackage{geometry}
\usepackage{latexsym}
\usepackage{mathrsfs}
\usepackage{pifont}
\usepackage{makecell}
\usepackage{amsfonts}
\usepackage{amssymb}
\usepackage{microtype}
\usepackage{subfigure}    %for parallel figures
\usepackage{graphicx}
\usepackage{epstopdf}
\usepackage{float}
\usepackage{multirow}
\usepackage{bm}
\usepackage{amsmath}
\usepackage{amsthm}
\usepackage{color}
\usepackage[subnum]{cases}
\usepackage{rotating}
\usepackage{enumitem}
\usepackage{accents}
\usepackage{algorithm}
\usepackage{algorithmic}
\usepackage[title]{appendix}
\usepackage{mathtools}
\usepackage{empheq}
\usepackage{caption}
\usepackage{url}
\usepackage{booktabs}
\usepackage{tabularx}
\usepackage{tcolorbox}
\usepackage{verbatim}
\usepackage[normalem]{ulem}

\makeatletter
\@addtoreset{equation}{section}
\makeatother

\newcommand{\email}[1]{\protect\href{mailto:#1}{#1}}

\definecolor{blue}{rgb}{0,0,0.9}
\definecolor{red}{rgb}{0.9,0,0}
\definecolor{green}{rgb}{0,0.9,0}
\definecolor{brown}{rgb}{0.6,0.1,0.1}
\definecolor{lightgreen}{rgb}{0.1,0.5,0.1}

\usepackage[colorlinks=true,
breaklinks=true,
bookmarks=true,
urlcolor=blue,
citecolor=lightgreen,
linkcolor=lightgreen,
bookmarksopen=false,
draft=false]{hyperref}

\graphicspath{{figures/}}

\title{A Variable-metric Inexact Progressive Hedging Algorithm \\ with a Lagged Stopping Criterion for \\ Multistage Stochastic Variational Inequalities}

\author{Min Zhang\thanks{School of Mathematics and Information Science, Guangzhou University (\email{zhangmin1206@gzhu.edu.cn})}
\and Lei Yang\thanks{(Corresponding author) School of Computer Science and Engineering, Sun Yat-sen University (\email{yanglei39@mail.sysu.edu.cn})}
\and Jie Sun\thanks{School of EECMS, Curtin University (\email{jie.sun@curtin.edu.au})}
}

\date{Updated on \today}

\begin{document}

\newtheorem{property}{Property}[section]
\newtheorem{proposition}{Proposition}[section]
\newtheorem{append}{Appendix}[section]
\newtheorem{definition}{Definition}[section]
\newtheorem{lemma}{Lemma}[section]
\newtheorem{corollary}{Corollary}[section]
\newtheorem{theorem}{Theorem}[section]
\newtheorem{remark}{Remark}[section]
\newtheorem{problem}{Problem}[section]
\newtheorem{example}{Example}[section]
\newtheorem{assumption}{Assumption}
\renewcommand*{\theassumption}{\Alph{assumption}}

\maketitle

\begin{abstract}
The progressive hedging algorithm (PHA) solves multistage stochastic variational inequalities (SVIs) by independently solving regularized scenario subproblems and enforcing nonanticipativity through primal-dual updates. Since these subproblems are often solved iteratively, practical implementations need stopping criteria that control errors while preserving convergence. We propose iVMPHA-Lag, a variable-metric inexact PHA with a lagged stopping criterion for monotone SVIs. The method allows scenario block metrics to vary across outer iterations while remaining compatible with nonanticipativity. For each scenario, an error budget combines a positive fraction of the previous budget with a progress measure from the preceding outer iteration. This gives a stopping tolerance fixed before the inner iterations begin, so each scenario can check its criterion independently without a prescribed summable tolerance sequence. With positive initial budgets and a positive inexactness parameter, we prove finite inner termination provided the inner solver produces valid approximations whose residual and enlargement errors tend to zero. Under suitable conditions, we establish global convergence of the primal-dual iterates and local $R$-linear convergence under metric subregularity at the limit point. Numerical experiments on synthetic strongly monotone box-constrained SVIs show that iVMPHA-Lag with adaptive metrics uses fewer scenario evaluations and less time on average than exact PHA and the competing inexact methods among the tested configurations in all four baseline regimes. Its average time advantage over the competing inexact methods persists across the tested problem sizes. Further comparisons show that metric adaptation improves efficiency over the fixed scalar configuration, while lagged stopping and positive budget retention reduce inner computational work.

\vspace{5mm}
\end{abstract}

%%%%%%%%%%%%%%%%%%%%%%%%%%%%%%%%%%%%%%%%%%%%%%%%%%%%%%%%%%%%
\section{Introduction}\label{sec:introduction}

Stochastic variational inequalities (SVIs) provide a unified framework for modeling multistage optimization and equilibrium problems under uncertainty, encompassing multistage stochastic programs and noncooperative games in which decision makers choose their actions sequentially in response to progressively revealed information \cite{rw2017svi,rs2019,zsx2019,sunchen2021,jiangsun2023}. This framework has a wide range of applications, including electricity market equilibria \cite{philpott2016equilibrium}, production and supply planning \cite{jiang2019regularized}, traffic assignment \cite{wang2023solving}, and medical supply management \cite{li2022two}.

We consider a finite scenario tree with $s$ scenarios. Let $\bm{x}_{i}\in\mathbb{R}^{d}$ denote the decision vector along scenario $i$, and write $\bm{x}=(\bm{x}_{1}^{\top},\ldots,\bm{x}_{s}^{\top})^{\top}$. Each scenario vector comprises decisions at successive stages, and decisions at the same stage must agree across scenarios that share the same information history up to that stage. These consistency requirements define the nonanticipativity subspace $\mathcal{N}$. Let $\mathcal{M}:=\mathcal{N}^{\perp}$, where orthogonality is taken with respect to the probability-weighted inner product. Under a standard constraint qualification, the multistage stochastic variational inequality (SVI) can be equivalently written in the following primal-dual extensive form:
\begin{equation}\label{SVI}
\text{Find} \quad \bm{x}\in\mathcal{N} \quad \text{and} \quad \bm{w}\in \mathcal{M} \quad\text{such that}\quad \bm{w}_{i}\in\mathcal{T}_{i}(\bm{x}_{i}), \quad \forall~i=1,\ldots,s.
\end{equation}
Here, the vector $\bm{w}=(\bm{w}_{1}^{\top},\ldots,\bm{w}_{s}^{\top})^{\top}$ represents the multipliers associated with the nonanticipativity constraints, and each scenario operator has the form $\mathcal{T}_{i}:=F_{i}+N_{C_{i}}$, where $F_{i}:\mathbb{R}^{d}\to\mathbb{R}^{d}$ is continuous and monotone, $C_{i}\subseteq\mathbb{R}^{d}$ is a nonempty closed convex feasible set, and $N_{C_{i}}$ is its Euclidean normal-cone operator. Section \ref{subsec:primal-dual} gives the detailed problem formulation and the constraint qualification under which the formulation is equivalent to \eqref{SVI}. Although the operators $\mathcal{T}_{i}$ act separately on the scenario decisions, the nonanticipativity constraint $\bm{x}\in\mathcal{N}$ couples these decisions across scenarios, making a direct solution of the full problem computationally demanding as the scenario tree grows.

The progressive hedging algorithm (PHA) exploits this structure by separating the scenario calculations while enforcing nonanticipativity through a coordination step \cite{rs2019,rw1991}. Given a regularization parameter $r>0$ and an iterate $(\bm{x}^{k},\bm{w}^{k})\in\mathcal{N}\times\mathcal{M}$, write $\widehat{\bm{x}}^{k+1} =(\widehat{\bm{x}}_{i}^{k+1})_{i=1}^{s}$. One exact PHA iteration consists of the following steps:
\begin{equation}\label{eq:classical-pha}
\left\{\begin{aligned} \textbf{Step 1: } &\text{Solve }~ \widehat{\bm{x}}_{i}^{k+1} \;\text{ via }\; 0\in {\mathcal{T}}_{i}(\widehat{\bm{x}}^{k+1}_{i}) -\bm{w}^{k}_{i} +r(\widehat{\bm{x}}^{k+1}_{i}-\bm{x}^{k}_{i}), \;\forall~i,  \\[0.3em]
\textbf{Step 2: } &\text{Update } \bm{x}^{k+1}=P_{\mathcal{N}}(\widehat{\bm{x}}^{k+1}) \text{ and } \bm{w}^{k+1} = \bm{w}^{k}-r(\widehat{\bm{x}}^{k+1}-\bm{x}^{k+1}),
\end{aligned}\right.
\end{equation}
where $P_{\mathcal{N}}$ denotes the probability-weighted orthogonal projection onto $\mathcal{N}$. The regularized scenario subproblems in Step~1 can be solved independently and in parallel. Step 2 then restores nonanticipativity by projection and updates the dual variable to account for disagreement among the scenario solutions.

The progressive decoupling algorithm (PDA) generalizes the PHA to linkage problems involving maximal monotone operators and linear coupling subspaces \cite{rockafellar2019progressive}. It includes \eqref{eq:classical-pha} as a special case, and subsequent developments allow variable metrics compatible with the coupling subspaces \cite{rockafellar2024generalizations}. Progressive hedging and decoupling have also been developed and analyzed in nonmonotone settings under elicitable monotonicity and more general assumptions \cite{rockafellar2019progressive-1,sun2021elicited,cui2024localized,evens2026spingarn}. In practice, the regularized subproblems often have no closed-form solutions and are therefore solved iteratively. Since requiring exact solutions at every outer iteration is generally impractical, it is computationally crucial to develop inexact versions of PHA/PDA, in which computable stopping criteria are used to control the errors in the approximate subproblem solutions while preserving convergence of the outer iterates.

For the scenario structure considered here, let $B_{i,k}\in\mathbb{R}^{d\times d}$ be a symmetric positive definite matrix that determines the proximal regularization for scenario $i$ at the $k$-th outer iteration. To describe the errors arising when solving a scenario subproblem, we consider the following relation
\begin{equation*}
\bm{\vartheta}_{i}^{k+1} \in \mathcal{T}_{i}^{\varepsilon_{i,k+1}}(\bm{y}_{i}^{k+1}) -\bm{w}_{i}^{k} + B_{i,k}^{-1}(\bm{y}_{i}^{k+1}-\bm{x}_{i}^{k}),
\end{equation*}
where $\bm{y}_{i}^{k+1}\in\mathbb{R}^{d}$ is the approximate scenario solution, $\bm{\vartheta}_{i}^{k+1}\in\mathbb{R}^{d}$ is the residual in the regularized inclusion with the enlarged operator, and $\varepsilon_{i,k+1}\geq0$ quantifies the relaxation of the operator relation through enlargement. The enlarged operator $\mathcal{T}_{i}^{\varepsilon_{i,k+1}}$ is defined in Section \ref{subsec:operators} and satisfies $\mathcal{T}_{i}^{0}=\mathcal{T}_{i}$. Thus, setting both $\bm{\vartheta}_{i}^{k+1}=0$ and $\varepsilon_{i,k+1}=0$ recovers the exact regularized subproblem. The choice $B_{i,k}=r^{-1}I$, where $I$ is the identity matrix, gives the scalar regularization used in \eqref{eq:classical-pha}.

A stopping criterion terminates the inner iterations once a computable error measure does not exceed the chosen tolerance. Existing criteria differ in how this tolerance is determined. We refer to criteria that require an externally prescribed summable tolerance sequence as \emph{absolute-type} criteria. Rockafellar's variable-metric inexact PDA \cite{rockafellar2026advances,rockafellar2024generalizations} uses two such criteria with zero enlargement error. The first directly bounds the residual with the prescribed tolerance, while the second scales the prescribed tolerance by the distance between the approximate solution and the current iterate, capped at one. Under suitable metric assumptions, the first criterion guarantees global convergence, while the second yields Q-linear convergence of the variable-metric distance to the solution set under metric subregularity. Li, Rockafellar, and Sun \cite{LiRockafellarSun2026} recently extended the analysis to semidefinite metrics while retaining summable tolerance requirements. However, choosing a summable sequence in advance can be difficult, since tight tolerances may cause unnecessary inner work, whereas loose tolerances may slow the outer iteration.

\emph{Relative-type} criteria use a fixed inexactness parameter to set the tolerances based on quantities computed during the algorithm, without prescribing a summable tolerance sequence. Solodov and Svaiter \cite{solodov1999hybrid} developed this approach in their hybrid proximal extragradient (HPE) framework; see \cite{buscaglia2025general,lotito2009class,monteiro2010complexity} for subsequent analyses and variable-metric extensions. For monotone SVIs, Buscaglia, Lotito, and Parente \cite{buscaglia2024inexact} proposed an inexact PHA with fixed scalar regularization and an aggregate relative-error criterion based on current scenario candidates. They established global convergence and, under outer Lipschitz continuity of the inverse of the associated scaled partial-inverse operator, linear convergence. By comparison, Burachik, Sagastiz\'{a}bal, and Scheimberg \cite{burachik2006inexact} developed an inexact partial-inverse method with a relative-type stopping criterion that can be checked separately for each subproblem. Through the scaled partial-inverse formulation underlying PHA, this method can be adapted to SVIs. Marques Alves and Costa Lima \cite{costa2017inexact} later clarified the admissible range of its inexactness parameter and established pointwise and ergodic iteration-complexity bounds for related inexact partial-inverse methods. However, the inexact PHA variants considered above use fixed scalar regularization, and their stopping tolerances depend on the current candidate quantities. This motivates the following question:

\begin{tcolorbox}[colback=black!0!white]
\emph{Can we design a variable-metric inexact PHA with relative-type error control, whose scenario stopping tolerances are fixed before the inner iterations begin and can be checked independently, without prescribing a summable tolerance sequence, while still guaranteeing convergence of the outer iterates?}
\end{tcolorbox}

We answer this question by proposing iVMPHA-Lag, a variable-metric inexact PHA with lagged scenario error budgets. The main contributions are summarized as follows.
\begin{enumerate}[label=(\arabic*),leftmargin=*]

\item We develop a lagged relative-type stopping criterion that controls both residual and enlargement errors separately for each scenario. Each error budget combines progress from the preceding outer iteration with a retained fraction of the previous budget, giving a stopping tolerance fixed throughout the inner iterations. The algorithm allows scenario block metrics compatible with the nonanticipativity subspace and uses a residual-corrected primal update.

\item We prove finite inner termination with positive initial budgets, positive retention, and a positive inexactness parameter, provided the subsolvers produce candidate triplets satisfying the scenario inclusions with residual and enlargement errors tending to zero. Under suitable conditions, we establish global primal-dual convergence and summability of the aggregate error budgets, squared residual norms, and enlargement errors. We further prove local $R$-linear convergence of the primal-dual iterates under metric subregularity at the limit point.

\item We compare iVMPHA-Lag with exact PHA and existing inexact PHA methods on synthetic multistage SVIs. Among the tested configurations, iVMPHA-Lag with adaptive metrics achieves the lowest mean scenario-evaluation count and elapsed time in all four baseline regimes and maintains an average time advantage over the competing inexact methods across the tested problem sizes. Further experiments show that the lagged stopping and positive budget retention reduce inner computational work.

\end{enumerate}

The rest of the paper is organized as follows. Section \ref{sec:model} introduces notation and preliminaries. Section \ref{sec:iPDA} presents iVMPHA-Lag and compares its stopping criterion with existing ones. Section \ref{sec:convergence} establishes global convergence and local $R$-linear convergence. Section \ref{sec:numerical} reports the numerical experiments, and Section \ref{sec:conclusion} concludes the paper. Appendix \ref{app:vm-hpe} relates iVMPHA-Lag to the variable-metric HPE framework when the enlargement error is zero.

%%%%%%%%%%%%%%%%%%%%%%%%%%%%%%%%%%%%%%%%%%%%%%%%%%%%%%%%%%%%
\section{Notation and preliminaries}\label{sec:model}

This section introduces the multistage SVI formulation, monotone operators and their enlargements, the partial-inverse representation, and the metric assumptions used in the analysis.
Scalars, vectors, and matrices are denoted by lowercase, bold lowercase, and uppercase letters, respectively. We use $\mathbb{R}$, $\mathbb{R}_{+}$, $\mathbb{R}^{n}$, and $\mathbb{R}^{m\times n}$ to denote the sets of real numbers, nonnegative real numbers, $n$-dimensional real vectors, and $m\times n$ real matrices, respectively. Unless otherwise specified, $\langle\cdot,\cdot\rangle$ and $\|\cdot\|$ are the Euclidean inner product and norm, and $I$ is the identity matrix of the appropriate size. The operator $\operatorname{Diag}$ forms a block diagonal matrix from a list of matrices or a diagonal matrix from a vector.

%%%%%%%%%%%%%%%%%%%%%%%%%%%%%%%%%%%%%%%%%%%%%%%%%%%%%%%%%%%%
\subsection{\texorpdfstring{Multistage SVI}{Multistage SVI}}\label{subsec:primal-dual}

\noindent Consider a finite scenario tree with $T$ stages and $s$ scenarios. The probability of scenario $i$ is $p_{i}>0$, with $\sum_{i=1}^{s}p_{i}=1$. The decision vector along scenario $i$ is $\bm{x}_{i}=(\bm{x}_{i,1}^{\top},\ldots,\bm{x}_{i,T}^{\top})^{\top}$, where $\bm{x}_{i,t}\in\mathbb{R}^{d_{t}}$ contains the decisions at stage $t$ and $d=\sum_{t=1}^{T}d_{t}$. The full decision vector $\bm{x}=(\bm{x}_{1}^{\top},\ldots,\bm{x}_{s}^{\top})^{\top}$ belongs to the decision space $\mathcal{H}=\prod_{i=1}^{s}\mathbb{R}^{d}$, equipped with the probability-weighted inner product
\begin{equation*}
\langle\bm{x},\bm{y}\rangle_{p}:=\sum_{i=1}^{s}p_{i}\langle\bm{x}_{i},\bm{y}_{i}\rangle,\quad \forall~\bm{x},\bm{y}\in\mathcal{H}.
\end{equation*}
The induced norm $\|\cdot\|_{p}$ is defined by $\|\bm{x}\|_{p}^{2}:=\langle\bm{x},\bm{x}\rangle_{p}=\sum_{i=1}^{s}p_{i}\|\bm{x}_{i}\|^{2}$. In this way, $\mathcal{H}$ forms a finite-dimensional Hilbert space. For primal-dual pairs in $\mathcal{H}\times\mathcal{H}$, we use the subscript $p$ for the product norm defined by $\|(\bm{x},\bm{w})\|_{p}^{2}:=\|\bm{x}\|_{p}^{2}+\|\bm{w}\|_{p}^{2}.$ For $\bm{v}\in\mathcal{H}$ and $\mathcal{S}\subseteq\mathcal{H}$, the distance from $\bm{v}$ to $\mathcal{S}$ is defined by $\operatorname{dist}_{p}(\bm{v},\mathcal{S}):=\inf_{\bm{z}\in\mathcal{S}}\|\bm{v}-\bm{z}\|_{p}$, with $\operatorname{dist}_{p}(\bm{v},\varnothing)=+\infty$. The same definition applies on $\mathcal{H}\times\mathcal{H}$ using the product norm.

Nonanticipativity requires stage-$t$ decisions to agree across scenarios that share the same information history when the stage-$t$ decision is made. Let $\mathcal{N}$ be the subspace of all decision vectors satisfying these requirements, and let $\mathcal{M}:=\mathcal{N}^{\perp}$ be its orthogonal complement with respect to $\langle\cdot,\cdot\rangle_{p}$. Then, $\mathcal{H}=\mathcal{N}\oplus\mathcal{M}$. We denote the corresponding orthogonal projections by $P_{\mathcal{N}}$ and $P_{\mathcal{M}}$, respectively. The projection $P_{\mathcal{N}}$ is computed by taking conditional expectations at each stage \cite{rs2019}. The complementary projection is $P_{\mathcal{M}}=I-P_{\mathcal{N}}$.

For each scenario $i$, let $C_{i}\subseteq\mathbb{R}^{d}$ be a nonempty closed convex feasible set, and let $F_{i}:\mathbb{R}^{d}\to\mathbb{R}^{d}$ be continuous and monotone. Define $\mathcal{C}:=\prod_{i=1}^{s}C_{i}$ and $\mathcal{F}(\bm{x}):=(F_{1}^{\top}(\bm{x}_{1}),\ldots,F_{s}^{\top}(\bm{x}_{s}))^{\top}$. The set $\mathcal{C}\cap\mathcal{N}$ consists of decisions that satisfy both the scenario constraints and nonanticipativity. In its basic variational form, the multistage stochastic VI is given by
\begin{equation}\label{eq:svi-basic}
\text{Find}\quad \bm{x}\in\mathcal{C}\cap\mathcal{N} \quad\text{such that}\quad \langle \mathcal{F}(\bm{x}),\bm{z}-\bm{x}\rangle_{p}\geq0 \quad\text{for every }\bm{z}\in\mathcal{C}\cap\mathcal{N}.
\end{equation}

To obtain a primal-dual representation of \eqref{eq:svi-basic}, define the scenario operator $\mathcal{T}_{i}:=F_{i}+N_{C_{i}}$. Here $N_{C_{i}}$ is the Euclidean normal-cone operator, defined at $\bm{y}_{i}\in C_{i}$ by $N_{C_{i}}(\bm{y}_{i}) = \left\{ \bm{v}_{i}\in\mathbb{R}^{d}: \langle\bm{v}_{i},\bm{z}_{i}-\bm{y}_{i}\rangle\leq0 \text{ for all }\bm{z}_{i}\in C_{i} \right\}$, and $N_{C_{i}}(\bm{y}_{i})=\varnothing$ otherwise. Define the product operator $\mathcal{T}(\bm{x}):=\prod_{i=1}^{s}\mathcal{T}_{i}(\bm{x}_{i})$. Under a constraint qualification such as $\mathcal{N}\cap\operatorname{ri}(\mathcal{C})\neq\varnothing$, where $\operatorname{ri}$ denotes the relative interior, \eqref{eq:svi-basic} is equivalent to the primal-dual formulation \eqref{SVI} \cite{rw2017svi,rs2019}, where $\bm{w}\in\mathcal{T}(\bm{x})$ collects the scenario inclusions $\bm{w}_{i}\in\mathcal{T}_{i}(\bm{x}_{i})$. We denote its solution set by
\begin{equation*}
\Omega := \left\{ (\bm{x},\bm{w})\in\mathcal{N}\times\mathcal{M}: \bm{w}\in\mathcal{T}(\bm{x}) \right\},
\end{equation*}
and assume throughout that $\Omega\neq\varnothing$.

%%%%%%%%%%%%%%%%%%%%%%%%%%%%%%%%%%%%%%%%%%%%%%%%%%%%%%%%%%%%
\subsection{Monotone operators and enlargements}\label{subsec:operators}

\noindent We next recall the definitions of monotone operators and their enlargements.

\begin{definition}
Let $\mathcal{A}:\mathcal{H}\rightrightarrows\mathcal{H}$ be a set-valued operator, with $\operatorname{gph}\mathcal{A} := \{(\bm{y},\bm{u}):\bm{u}\in\mathcal{A}(\bm{y})\}.$ The operator $\mathcal{A}$ is monotone if $\langle \bm{y}-\bm{z},\bm{u}-\bm{v}\rangle_{p} \geq 0$ for any $(\bm{y},\bm{u}),~(\bm{z},\bm{v})\in\operatorname{gph}\mathcal{A}$.
%\begin{equation*}
%\langle \bm{y}-\bm{z},\bm{u}-\bm{v}\rangle_{p} \geq 0,\quad \forall \; (\bm{y},\bm{u}),~(\bm{z},\bm{v})\in\operatorname{gph}\mathcal{A}.
%\end{equation*}
It is maximally monotone if it is monotone and its graph is not properly contained in the graph of any other monotone operator.
\end{definition}

Since $F_{i}$ is continuous and monotone on $\mathbb{R}^{d}$ and $N_{C_{i}}$ is maximally monotone, the standard sum theorem implies that $\mathcal{T}_{i}=F_{i}+N_{C_{i}}$ is maximally monotone. Consequently, the product operator $\mathcal{T}$ is maximally monotone on the probability-weighted space $\mathcal{H}$.

\begin{definition}\cite{solodov1999hybrid}
Let $\mathcal{A}$ be a maximally monotone operator on a Hilbert space with inner product $\langle\cdot,\cdot\rangle$. For $\varepsilon\geq0$, its \emph{$\varepsilon$-enlargement} is defined by
\begin{equation*}
\mathcal{A}^{\varepsilon}(\bm{y}) := \left\{ \bm{u}: \langle \bm{y}-\bm{z},\bm{u}-\bm{v}\rangle \geq-\varepsilon \quad \text{for all }(\bm{z},\bm{v})\in\operatorname{gph}\mathcal{A} \right\}.
\end{equation*}
\end{definition}

By maximal monotonicity, $\mathcal{A}^{0}=\mathcal{A}$, so setting $\varepsilon=0$ recovers exact membership in the graph of $\mathcal{A}$. For each scenario operator $\mathcal{T}_{i}$, the enlargement is defined using the Euclidean inner product on $\mathbb{R}^{d}$. The enlargement of the product operator $\mathcal{T}$ instead uses the probability-weighted inner product $\langle\cdot,\cdot\rangle_{p}$. Moreover, suppose that each scenario provides a pair $(\bm{y}_{i},\bm{u}_{i})$ and an error $\varepsilon_{i}\geq0$ satisfying $\bm{u}_{i}\in\mathcal{T}_{i}^{\varepsilon_{i}}(\bm{y}_{i})$, for all $i=1,\ldots,s$. Set $\bm{y}=(\bm{y}_{i})_{i=1}^{s}$ and $\bm{u}=(\bm{u}_{i})_{i=1}^{s}$. Multiplying the scenario enlargement inequalities by $p_{i}$ and summing over $i$ gives
\begin{equation}\label{eq:weighted-enlargement}
\bm{u}\in\mathcal{T}^{\varepsilon}(\bm{y}) \quad \text{with} \quad \varepsilon:=\sum_{i=1}^{s} p_{i}\varepsilon_{i}.
\end{equation}
Relation \eqref{eq:weighted-enlargement} combines the scenario inclusions into an inclusion for $\mathcal{T}$, with enlargement error equal to the probability-weighted sum of the scenario enlargement errors.

%%%%%%%%%%%%%%%%%%%%%%%%%%%%%%%%%%%%%%%%%%%%%%%%%%%%%%%%%%%%
\subsection{\texorpdfstring{Partial inverse and the exact variable-metric PDA}{Partial inverse and exact variable-metric PDA}}\label{subsec:partial-inverse}

\noindent Spingarn's partial inverse \cite{spingarn1983partial} allows us to rewrite the primal-dual formulation \eqref{SVI} of the multistage SVI as a zero-finding problem in a single variable. For the maximally monotone operator $\mathcal{T}$ and the subspace $\mathcal{N}$, the partial inverse $\mathcal{T}_{\mathcal{N}}:\mathcal{H}\rightrightarrows\mathcal{H}$ is defined, for $\bm{z},\bm{v}\in\mathcal{H}$, by
\begin{equation}\label{eq:partial-inverse}
\bm{v}\in\mathcal{T}_{\mathcal{N}}(\bm{z}) \quad\Longleftrightarrow\quad P_{\mathcal{N}}\bm{v}+P_{\mathcal{M}}\bm{z} \in \mathcal{T}\!\left( P_{\mathcal{N}}\bm{z}+P_{\mathcal{M}}\bm{v} \right).
\end{equation}
The operator $\mathcal{T}_{\mathcal{N}}$ is maximally monotone. Its enlargement is obtained by applying the same partial-inverse construction to $\mathcal{T}^{\varepsilon}$ \cite{costa2017inexact,spingarn1983partial}:
\begin{equation}\label{eq:partial-enlargement}
(\mathcal{T}_{\mathcal{N}})^{\varepsilon} = (\mathcal{T}^{\varepsilon})_{\mathcal{N}}.
\end{equation}

The orthogonal decomposition $\mathcal{H}=\mathcal{N}\oplus\mathcal{M}$ gives the representation $\bm{z}=\bm{x}+\bm{w}$, with $\bm{x}=P_{\mathcal{N}}\bm{z}$ and $\bm{w}=P_{\mathcal{M}}\bm{z}$. Since $\|\bm{z}\|_{p}^{2}=\|\bm{x}\|_{p}^{2}+\|\bm{w}\|_{p}^{2}$, this correspondence is an isometry from $\mathcal{N}\times\mathcal{M}$ onto $\mathcal{H}$. Setting $\bm{v}=0$ in \eqref{eq:partial-inverse} yields
\begin{equation*}
(\bm{x},\bm{w})\in\Omega \quad\Longleftrightarrow\quad 0\in\mathcal{T}_{\mathcal{N}}(\bm{x}+\bm{w}).
\end{equation*}
Therefore, solving the primal-dual formulation \eqref{SVI} of the multistage SVI is equivalent to finding a zero of the maximally monotone operator $\mathcal{T}_{\mathcal{N}}$. This allows us to apply the proximal point algorithm (PPA) for solving SVIs.

To retain independent scenario subproblems, let $B_{k}=\operatorname{Diag}(B_{1,k},\ldots,B_{s,k})$ have symmetric positive definite blocks $B_{i,k}\in\mathbb{R}^{d\times d}$ and satisfy $B_{k}\mathcal{N}\subseteq\mathcal{N}$. We call such metrics \emph{compatible} with $\mathcal{N}$. By self-adjointness, $B_{k}$ also preserves $\mathcal{M}$, so $B_{k}$ and $B_{k}^{-1}$ commute with both projections, and the operator
\begin{equation*}
M_{k} := B_{k}^{-1}P_{\mathcal{N}} + B_{k}P_{\mathcal{M}}
\end{equation*}
is positive definite and self-adjoint. With $\bm{z}^{k}=\bm{x}^{k}+\bm{w}^{k}$, the corresponding variable-metric PPA iteration is
\begin{equation}\label{eq:exact-ppa}
0 \in \mathcal{T}_{\mathcal{N}}(\bm{z}^{k+1}) + M_{k}(\bm{z}^{k+1}-\bm{z}^{k}).
\end{equation}

Using the definition of the partial inverse and the commutation of $B_{k}$ with the projections, we can express \eqref{eq:exact-ppa} as the exact variable-metric PDA below, with $\widehat{\bm{x}}^{k+1}=(\widehat{\bm{x}}_{i}^{k+1})_{i=1}^{s}$:
\begin{equation}\label{eq:block-pda}
\left\{\begin{aligned} \textbf{Step 1: } &\text{Solve }~ \widehat{\bm{x}}_{i}^{k+1} \;\text{ via }\; 0\in {\mathcal{T}}_{i}(\widehat{\bm{x}}^{k+1}_{i}) -\bm{w}^{k}_{i} + B_{i,k}^{-1}(\widehat{\bm{x}}^{k+1}_{i}-\bm{x}^{k}_{i}), \;\forall~i,  \\[0.3em]
\textbf{Step 2: } &\text{Update } \bm{x}^{k+1}=P_{\mathcal{N}}(\widehat{\bm{x}}^{k+1}),\; \bm{w}^{k+1} = \bm{w}^{k}-B_{k}^{-1}(\widehat{\bm{x}}^{k+1}-\bm{x}^{k+1}).
\end{aligned}\right.
\end{equation}
Thus, exact variable-metric PDA is precisely a variable-metric PPA applied to Spingarn's partial inverse $\mathcal{T}_{\mathcal{N}}$, expressed in the original primal-dual variables \cite{sun2021elicited}.

%%%%%%%%%%%%%%%%%%%%%%%%%%%%%%%%%%%%%%%%%%%%%%%%%%%%%%%%%%%%
\subsection{\texorpdfstring{Metric assumptions and norm relations}{Metric assumptions and norm relations}}\label{subsec:metrics}

\noindent The preceding formulation uses positive definite scenario blocks and compatibility with $\mathcal{N}$. For the convergence analysis, we further require the eigenvalues of $B_{k}$ to be uniformly bounded away from zero and infinity, and control how much the metric changes between successive iterations. For self-adjoint operators $H_{1}$ and $H_{2}$ on $\mathcal{H}$, we write $H_{1}\preceq H_{2}$ if $H_{2}-H_{1}$ is positive semidefinite.

\begin{assumption}[Scenario-compatible metrics]\label{assmp:Bk}
For every $k$, let $B_{i,k}\in\mathbb{R}^{d\times d}$ be symmetric positive definite matrices. Set $B_{k}=\operatorname{Diag}(B_{1,k},\ldots,B_{s,k})$ and assume that $B_{k}\mathcal{N}\subseteq\mathcal{N}$. There exist constants $0<\underline{\lambda}\leq\overline{\lambda}<\infty$ and a sequence $\{\eta_{k}\}\subseteq\mathbb{R}_{+}$ such that
\begin{equation*}
\underline{\lambda}I \preceq B_{k} \preceq\overline{\lambda}I, \quad \frac{1}{1+\eta_{k}}B_{k} \preceq B_{k+1} \preceq (1+\eta_{k})B_{k}, \quad \sum_{k=0}^{\infty}\eta_{k}<\infty.
\end{equation*}
\end{assumption}

A simple compatible choice is to use the same metric in every scenario, with one block for each stage:
\begin{equation}\label{stage-metric}
B_{i,k}=\operatorname{Diag}(b_{1,k} I_{d_{1}}, \ldots, b_{T,k} I_{d_{T}}), \quad b_{t,k}>0,\quad i=1,\ldots,s.
\end{equation}
Here, to specify the dimension, $I_{d_{t}}$ denotes the identity matrix in $\mathbb{R}^{d_{t}\times d_{t}}$. For this choice, decisions that agree at a stage remain equal after multiplication by $B_{k}$, so $B_{k}\mathcal{N}\subseteq\mathcal{N}$. More generally, the block at stage $t$ may depend on the history before that stage, provided scenarios sharing that history use the same block. This gives a compatible metric with different scalings at different nodes of the scenario tree.

We now introduce the norms used to measure errors in these metrics. For a positive definite self-adjoint operator $H$ on $\mathcal{H}$, let $\|\bm{v}\|_{H,p}^{2}:=\langle\bm{v},H\bm{v}\rangle_{p}.$ For a positive definite self-adjoint scenario block $H_{i}$ on $\mathbb{R}^{d}$, we similarly write $\|\bm{v}_{i}\|_{H_{i}}^{2}:=\langle\bm{v}_{i},H_{i}\bm{v}_{i}\rangle$. In the convergence analysis, the primal and dual errors are measured using $B_{k}^{-1}$ and $B_{k}$, respectively. We combine these two metrics in the product-space operator
\begin{equation*}
G_{k} :=
\begin{pmatrix}
B_{k}^{-1} & 0\\
0 & B_{k}
\end{pmatrix},
\end{equation*}
with the associated norm
\begin{equation}\label{eq:metric-norms}
\|(\bm{x},\bm{w})\|_{G_{k},p}^{2} := \|\bm{x}\|_{B_{k}^{-1},p}^{2} + \|\bm{w}\|_{B_{k},p}^{2}.
\end{equation}
We write $\operatorname{dist}_{G_{k},p}$ for the distance induced by $\|\cdot\|_{G_{k},p}$.

The spectral bounds in Assumption \ref{assmp:Bk} make $\|\cdot\|_{G_{k},p}$ uniformly equivalent to $\|\cdot\|_{p}$. The variation bound allows us to compare the same primal-dual error under two consecutive metrics. Since $\sum_{k=0}^{\infty}\eta_{k}<\infty$, the accumulated factors satisfy $\prod_{k=0}^{\infty}(1+\eta_{k})<\infty$. The following proposition records the norm relations used in the convergence analysis. The proof follows directly from Assumption \ref{assmp:Bk}, the definitions of $M_{k}$ and $G_{k}$, and the fact that $B_{k}$ preserves the orthogonal subspaces $\mathcal{N}$ and $\mathcal{M}$.

\begin{proposition}\label{prop:metric-identities}%[Metric relations and norm equivalence]
Suppose Assumption \ref{assmp:Bk} holds, and let $\Lambda:=\max\{\overline{\lambda},\underline{\lambda}^{-1}\}.$ Then, for every $k\geq0$, the following properties hold.
\begin{itemize}
\item[(i)] For every $(\bm{x},\bm{w})\in\mathcal{H}\times\mathcal{H}$,
\begin{equation}\label{eq:norm-equivalence}
\Lambda^{-1}\|(\bm{x},\bm{w})\|_{p}^{2} \leq \|(\bm{x},\bm{w})\|_{G_{k},p}^{2} \leq \Lambda\|(\bm{x},\bm{w})\|_{p}^{2},
\end{equation}
and
\begin{equation}\label{eq:successive-metric}
\|(\bm{x},\bm{w})\|_{G_{k+1},p}^{2} \leq (1+\eta_{k}) \|(\bm{x},\bm{w})\|_{G_{k},p}^{2}.
\end{equation}

\item[(ii)] For every $(\bm{x},\bm{w})\in\mathcal{N}\times\mathcal{M}$,
\begin{equation*}
\|(\bm{x},\bm{w})\|_{G_{k},p}^{2} = \|\bm{x}+\bm{w}\|_{M_{k},p}^{2}.
\end{equation*}

\item[(iii)] For $H\in\{B_{k},B_{k}^{-1}\}$ and every $\bm{v}\in\mathcal{H}$,
\begin{equation*}
\|H\bm{v}\|_{H^{-1},p} = \|\bm{v}\|_{H,p}, \qquad \|\bm{v}\|_{H,p}^{2} = \|P_{\mathcal{N}}\bm{v}\|_{H,p}^{2} + \|P_{\mathcal{M}}\bm{v}\|_{H,p}^{2}.
\end{equation*}
\end{itemize}
Consequently, for every nonempty set $\mathcal{S}\subseteq\mathcal{H}\times\mathcal{H}$,
\begin{equation}\label{eq:distance-equivalence}
\Lambda^{-1/2} \operatorname{dist}_{p}((\bm{x},\bm{w}),\mathcal{S}) \leq \operatorname{dist}_{G_{k},p}((\bm{x},\bm{w}),\mathcal{S}) \leq \Lambda^{1/2} \operatorname{dist}_{p}((\bm{x},\bm{w}),\mathcal{S}).
\end{equation}
\end{proposition}

%%%%%%%%%%%%%%%%%%%%%%%%%%%%%%%%%%%%%%%%%%%%%%%%%%%%%%%%%%%%
\section{An inexact progressive hedging algorithm with a lagged stopping criterion}\label{sec:iPDA}

\noindent In this section, we present iVMPHA-Lag, an inexact version of the variable-metric PHA \eqref{eq:block-pda} with a lagged stopping criterion. Each scenario uses the error budget updated at the preceding outer iteration to set its stopping tolerance for the current outer iteration. The complete iterative framework is given in Algorithm \ref{algo:iPDA}.

\begin{algorithm}[H]
\caption{iVMPHA-Lag for solving the multistage SVI \eqref{SVI}}\label{algo:iPDA}
Choose $(\bm{x}^{0},\bm{w}^{0})\in\mathcal{N}\times\mathcal{M}$ and $\sigma_{i,0}>0$ for $i=1,\ldots,s$. Fix an inexactness parameter $\rho\in[0,1)$ and a retention parameter $\delta\in(0,1-\rho^{2})$. Use metrics satisfying Assumption \ref{assmp:Bk}. For $k=0,1,\ldots$, perform the following steps:\vspace{2mm}

\textbf{Step 1. Inexact solution of scenario subproblems.} For each $i=1,\ldots,s$, find a triplet $({\bm{y}}_{i}^{k+1},{\bm{\vartheta}}_{i}^{k+1},\varepsilon_{i,k+1})\in \mathbb{R}^{d}\times\mathbb{R}^{d}\times\mathbb{R}_{+}$ satisfying
\begin{subnumcases}{\label{eq:block-inexcond}}
\bm{\vartheta}_{i}^{k+1} \in \mathcal{T}_{i}^{\varepsilon_{i,k+1}} (\bm{y}_{i}^{k+1}) -\bm{w}_{i}^{k} + B_{i,k}^{-1} (\bm{y}_{i}^{k+1}-\bm{x}_{i}^{k}),\label{eq:block-inexcond-a}\\[1mm]
\|\bm{\vartheta}_{i}^{k+1}\|_{B_{i,k}}^{2} +2\varepsilon_{i,k+1} \leq\rho^{2}\sigma_{i,k}^{2}.\label{eq:block-inexcond-b}
\end{subnumcases}
% \begin{subequations}
% \label{eq:block-inexcond}
% \begin{empheq}[left=\empheqlbrace]{align}
% &\quad \bm{\vartheta}_{i}^{k+1}\in\mathcal{T}_{i}^{\varepsilon_{i,k+1}}(\bm{y}_{i}^{k+1})-\bm{w}_{i}^{k}+B_{k,i}^{-1}(\bm{y}_{i}^{k+1}-\bm{x}_{i}^{k}),
% \label{eq:block-inexcond-a}\\[2mm]
% &\quad \|\bm{\vartheta}_{i}^{k+1}\|_{B_{k,i}}^{2}+2\varepsilon_{i,k+1}\leq\rho^{2}\sigma_{i,k}^{2}.
% \label{eq:block-inexcond-b}
% \end{empheq}
% \end{subequations}

\textbf{Step 2. Primal and dual updates.} Collect the accepted scenario vectors as $\bm{y}^{k+1}=(\bm{y}_{i}^{k+1})_{i=1}^{s}$ and $\bm{\vartheta}^{k+1}=(\bm{\vartheta}_{i}^{k+1})_{i=1}^{s}$. Compute $\bar{\bm{y}}^{k+1}=P_{\mathcal{N}}\bm{y}^{k+1}$, $\bar{\bm{\vartheta}}^{k+1}=P_{\mathcal{N}}\bm{\vartheta}^{k+1}$, and
\begin{equation}\label{eq:iPDA-update}
\bm{x}^{k+1} = \bar{\bm{y}}^{k+1}-B_{k}\bar{\bm{\vartheta}}^{k+1}, \quad \bm{w}^{k+1} = \bm{w}^{k}-B_{k}^{-1}(\bm{y}^{k+1}-\bar{\bm{y}}^{k+1}).
\end{equation}

\textbf{Step 3. Update the scenario error budgets.} For each $i=1,\ldots,s$, compute the progress measure $s_{i,k+1}^{2}$ and the next budget:
\begin{equation}\label{eq:budget-update}
\begin{aligned}
s_{i,k+1}^{2}&=\|\bar{\bm{y}}_{i}^{k+1}-\bm{x}_{i}^{k}\|_{B_{i,k}^{-1}}^{2}+ \|\bm{\vartheta}_{i}^{k+1}-\bar{\bm{\vartheta}}_{i}^{k+1}-B_{i,k}^{-1}(\bm{y}_{i}^{k+1}-\bar{\bm{y}}_{i}^{k+1})\|_{B_{i,k}}^{2},\\
\sigma_{i,k+1}^{2}&=s_{i,k+1}^{2}+\delta\sigma_{i,k}^{2}.
\end{aligned}
\end{equation}
\end{algorithm}

One can see from Algorithm \ref{algo:iPDA} that, in \textbf{Step 1}, each scenario $i$ uses the \emph{error budget} $\sigma_{i,k}^{2}$ from the preceding outer iteration to set its \emph{stopping tolerance} $\rho^{2}\sigma_{i,k}^{2}$. With this tolerance fixed, the inner solver computes a candidate triplet $(\bm{y}_{i}^{k+1},\bm{\vartheta}_{i}^{k+1},\varepsilon_{i,k+1})$ satisfying \eqref{eq:block-inexcond-a} and accepts it when the computable measure $\|\bm{\vartheta}_{i}^{k+1}\|_{B_{i,k}}^{2}+2\varepsilon_{i,k+1}$ of the residual and enlargement errors does not exceed this tolerance, as required by \eqref{eq:block-inexcond-b}. Once all scenario solutions have been accepted, \textbf{Step 2} projects $\bm{y}^{k+1}$ and $\bm{\vartheta}^{k+1}$ onto $\mathcal{N}$ to obtain $\bar{\bm{y}}^{k+1}$ and $\bar{\bm{\vartheta}}^{k+1}$, respectively. These projections are computed stage by stage by taking conditional expectations given the information available when the corresponding stage decision is made. The point $\bar{\bm{y}}^{k+1}$ is then corrected using the projected residual $\bar{\bm{\vartheta}}^{k+1}$ to obtain the new primal iterate $\bm{x}^{k+1}$, while the dual variable is updated using the deviations $\bm{y}^{k+1}-\bar{\bm{y}}^{k+1}$. Finally, \textbf{Step 3} computes the progress measure $s_{i,k+1}^{2}$ in \eqref{eq:budget-update} from the current scenario solutions and residuals, and adds the retained budget $\delta\sigma_{i,k}^{2}$ to form $\sigma_{i,k+1}^{2}$. This updated budget determines the stopping tolerance $\rho^{2}\sigma_{i,k+1}^{2}$ for the next outer iteration.

The following remark explains the roles of the metric, the retained budget, and the residual correction in these steps.

\begin{remark}\label{rem:zero-retention}
\begin{enumerate}[label=(\roman*)]
\item Since $B_{k}$ is block diagonal, the scenario subproblems in \textbf{Step 1} can be solved independently. Compatibility with $\mathcal{N}$ also gives $B_{k}\mathcal{N}\subseteq\mathcal{N}$ and $B_{k}^{-1}\mathcal{M}\subseteq\mathcal{M}$. Thus, in \textbf{Step 2}, the residual correction lies in $\mathcal{N}$ and the dual increment lies in $\mathcal{M}$, so $(\bm{x}^{k+1},\bm{w}^{k+1})\in\mathcal{N}\times\mathcal{M}$.

\item Retaining $\delta\sigma_{i,k}^{2}$ in \textbf{Step 3} ensures that $\sigma_{i,k}^{2}\geq\delta^{k}\sigma_{i,0}^{2}>0$. Hence, when $\rho>0$, each scenario has a positive stopping tolerance at every outer iteration. If $\delta=0$, the convergence estimates still hold for accepted sequences, but the next budget may vanish even when the current iterate is not a solution. For example, consider two equally likely scenarios with a shared scalar decision. Then $\mathcal{H}=\mathbb{R}^{2}$, $\mathcal{N}=\{(a,a)^{\top}:a\in\mathbb{R}\}$, and $\mathcal{M}=\{(a,-a)^{\top}:a\in\mathbb{R}\}$. Take $\mathcal{T}_{i}(\bm{x}_{i})=2 \bm{x}_{i}$, $B_{k}=I$, $\rho=1/2$, and initialize $\bm{x}^{0}=(1,1)^{\top}$, $\bm{w}^{0}=(0,0)^{\top}$, and $\sigma_{i,0}=4$. The candidate $\bm{y}^{1}=(1,1)^{\top}$, $\bm{\vartheta}^{1}=(2,2)^{\top}$ and $\varepsilon_{i,1}=0$ satisfies \eqref{eq:block-inexcond-a} and \eqref{eq:block-inexcond-b}. It yields $\bm{x}^{1}=(-1,-1)^{\top}$, $\bm{w}^{1}=(0,0)^{\top}$ and $s_{i,1}=0$ for both scenarios, although the unique primal-dual solution is $\bm{x}^{*}=\bm{0}$ and $\bm{w}^{*}=\bm{0}$. Without retention, the update $\sigma_{i,1}^{2}=s_{i,1}^{2}=0$ would require $\bm{\vartheta}_{i}^{2}=\bm{0}$ and $\varepsilon_{i,2}=0$, forcing the subproblem to be solved exactly at the next outer iteration. With positive retention, Algorithm \ref{algo:iPDA} instead gives $\sigma_{i,1}^{2}=16\delta>0$. The numerical experiments in Section \ref{subsec:budget-influence} also illustrate how positive budget retention can reduce the work required by the inner solvers.

\item With the residual correction $-B_{k}P_{\mathcal{N}}\bm{\vartheta}^{k+1}$ in \textbf{Step 2}, the primal-dual update satisfies the descent estimate in Lemma \ref{lem:one-step}. When $\rho=0$, the stopping criterion forces $\bm{\vartheta}_{i}^{k+1}=\bm{0}$ and $\varepsilon_{i,k+1}=0$ for every scenario. The correction then vanishes, and Algorithm \ref{algo:iPDA} reduces to the exact variable-metric PDA \eqref{eq:block-pda}. Taking $B_{k}=r^{-1}I$ further gives PHA \eqref{eq:classical-pha}.
\end{enumerate}
\end{remark}

We next describe how to compute the residual and enlargement error for an approximate scenario solution in \textbf{Step 1}. For $\mathcal{T}_{i}=F_{i}+N_{C_{i}}$, suppose that, at outer iteration $k$, the subsolver produces an approximate solution $\bm{y}_{i}^{k+1}\in C_{i}$ together with a normal vector $\bm{v}_{i}^{k+1}\in N_{C_{i}}(\bm{y}_{i}^{k+1})$. Since $F_{i}(\bm{y}_{i}^{k+1})+\bm{v}_{i}^{k+1}\in\mathcal{T}_{i}(\bm{y}_{i}^{k+1})$, we can set $\varepsilon_{i,k+1}=0$ and compute the residual
\begin{equation}\label{eq:normal-certificate}
\bm{\vartheta}_{i}^{k+1}=F_{i}(\bm{y}_{i}^{k+1})+\bm{v}_{i}^{k+1}-\bm{w}_{i}^{k}+B_{i,k}^{-1}(\bm{y}_{i}^{k+1}-\bm{x}_{i}^{k}).
\end{equation}
The resulting candidate triplet $(\bm{y}_{i}^{k+1},\bm{\vartheta}_{i}^{k+1},0)$ satisfies \eqref{eq:block-inexcond-a}. The stopping criterion \eqref{eq:block-inexcond-b} then compares the squared residual norm with the fixed tolerance $\rho^{2}\sigma_{i,k}^{2}$. In the numerical experiments of Section \ref{subsec:method-comparison}, the projection step of the diagonally preconditioned extragradient method provides $\bm{v}_{i}^{k+1}\in N_{C_{i}}(\bm{y}_{i}^{k+1})$. We then compute the residual by \eqref{eq:normal-certificate} and set $\varepsilon_{i,k+1}=0$.

When $C_{i}=K$ is a nonempty closed convex cone, we can choose the enlargement error using a complementarity gap. Let $K^{*}:=\{\bm{v}:\langle\bm{v},\bm{z}\rangle\geq0\text{ for all }\bm{z}\in K\}$ denote its dual cone. Suppose that the subsolver produces an approximate solution $\bm{y}_{i}^{k+1}\in K$ together with a vector $\bm{v}_{i}^{k+1}\in-K^{*}$. We take their complementarity gap as the enlargement error: $\varepsilon_{i,k+1}=-\langle\bm{y}_{i}^{k+1},\bm{v}_{i}^{k+1}\rangle\geq0.$ This choice gives $\bm{v}_{i}^{k+1}\in N_{K}^{\varepsilon_{i,k+1}}(\bm{y}_{i}^{k+1})$. By monotonicity of $F_{i}$, we have $F_{i}(\bm{y}_{i}^{k+1})+\bm{v}_{i}^{k+1}\in\mathcal{T}_{i}^{\varepsilon_{i,k+1}}(\bm{y}_{i}^{k+1})$. We then compute the residual by \eqref{eq:normal-certificate}, and the resulting candidate triplet satisfies \eqref{eq:block-inexcond-a}.

The following proposition establishes finite termination of \textbf{Step 1} when the candidate errors tend to zero.

\begin{proposition}[Finite termination of the inner iterations]\label{prop:finite-certification}
In Algorithm \ref{algo:iPDA}, let $0<\rho<1$ and fix an outer iteration $k$ and a scenario $i$. Suppose that the subsolver generates candidate triplets $(\bm{y}_{i}^{k+1,\ell},\bm{\vartheta}_{i}^{k+1,\ell},\varepsilon_{i,k+1,\ell})$, indexed by the inner iteration $\ell$, each satisfying \eqref{eq:block-inexcond-a} with $\varepsilon_{i,k+1,\ell}\geq0$. If
\begin{equation*}
\|\bm{\vartheta}_{i}^{k+1,\ell}\|_{B_{i,k}}^{2}+2\varepsilon_{i,k+1,\ell}\to0 \quad\text{as}\quad \ell\to\infty,
\end{equation*}
then, after finitely many inner iterations, a candidate triplet satisfies the stopping criterion \eqref{eq:block-inexcond-b} and is accepted. If every scenario subsolver satisfies these conditions at each outer iteration, \textbf{Step 1} terminates after finitely many inner iterations at every outer iteration.
\end{proposition}
\begin{proof}
The budget recursion \eqref{eq:budget-update}, together with $\sigma_{i,0}^{2}>0$ and $\delta>0$, gives $\sigma_{i,k}^{2}\geq\delta^{k}\sigma_{i,0}^{2}>0$. Since $\rho>0$, the stopping tolerance $\rho^{2}\sigma_{i,k}^{2}$ is positive and remains fixed throughout the inner iterations. As the error measure of the candidate sequence tends to zero, there exists an inner iteration $\ell_{i,k}$ such that
\begin{equation*}
\|\bm{\vartheta}_{i}^{k+1,\ell_{i,k}}\|_{B_{i,k}}^{2}+2\varepsilon_{i,k+1,\ell_{i,k}}\leq\rho^{2}\sigma_{i,k}^{2}.
\end{equation*}
The candidate triplet with index $\ell_{i,k}$ also satisfies \eqref{eq:block-inexcond-a} by assumption. A candidate triplet is therefore accepted no later than inner iteration $\ell_{i,k}$, and the inner solver terminates. Applying this argument to all $s$ scenarios proves finite termination of \textbf{Step 1} at each outer iteration under the stated conditions.
\end{proof}

%%%%%%%%%%%%%%%%%%%%%%%%%%%%%%%%%%%%%%%%%%%%%%%%%%%%%%%%%%%%
\subsection{Comparison with representative inexact PHA schemes}\label{subsec:relation-existing}

\noindent We next compare iVMPHA-Lag with representative absolute-type and relative-type PDA/PHA methods. For each method, we describe the candidate quantities, the stopping criterion used to accept them, and the resulting primal-dual update.

\begin{enumerate}[label=(\roman*),leftmargin=*]

\item \textit{iVMPDA-absI and iVMPDA-absII (absolute-type).} Rockafellar's variable-metric PDA \cite[Equations (2.10) and (2.22)]{rockafellar2026advances} uses two absolute-type stopping criteria. Each scenario subsolver provides a candidate pair $(\bm{y}_{i}^{k+1},\bm{\vartheta}_{i}^{k+1})$ satisfying \eqref{eq:block-inexcond-a} with $\varepsilon_{i,k+1}=0$. The aggregate residual must satisfy, respectively,
\begin{equation*}
\|\bm{\vartheta}^{k+1}\|_{B_{k},p} \leq a_{k}, \quad\text{or}\quad \|\bm{\vartheta}^{k+1}\|_{B_{k},p} \leq c_{k} \min\!\left\{ 1, \|\bm{y}^{k+1}-\bm{x}^{k}\|_{B_{k}^{-1},p} \right\},
\end{equation*}
where $\{a_{k}\}$ and $\{c_{k}\}$ are prescribed summable nonnegative sequences. The first bound is fixed before the inner iterations, whereas the second depends on the current candidate displacement and may change during them.

Both aggregate criteria admit scenario-wise sufficient conditions. For iVMPDA-absI, choose $a_{i,k}\geq 0$ with $\sum_{i=1}^{s}p_{i}a_{i,k}^{2}\leq a_{k}^{2}$ and require $\|\bm{\vartheta}_{i}^{k+1}\|_{B_{i,k}}\leq a_{i,k}$. For iVMPDA-absII, it suffices to require
\begin{equation*}
\|\bm{\vartheta}_{i}^{k+1}\|_{B_{i,k}} \leq c_{k}\min\{1,\|\bm{y}_{i}^{k+1}-\bm{x}_{i}^{k}\|_{B_{i,k}^{-1}}\}, \qquad i=1,\ldots,s.
\end{equation*}
Indeed, writing $d_{i}=\|\bm{y}_{i}^{k+1}-\bm{x}_{i}^{k}\|_{B_{i,k}^{-1}}$, we have
\begin{equation*}
\|\bm{\vartheta}^{k+1}\|_{B_{k},p}^{2} \leq c_{k}^{2}\sum_{i=1}^{s} p_{i}\min\{1,d_{i}^{2}\} \leq c_{k}^{2}\min\Bigl\{1,\sum_{i=1}^{s} p_{i} d_{i}^{2}\Bigr\}.
\end{equation*}

Once the scenario candidates have been accepted, both variants update
\begin{equation*}
\bm{x}^{k+1}=P_{\mathcal{N}}(\bm{y}^{k+1}) \text{ and } \bm{w}^{k+1} = \bm{w}^{k}-B_{k}^{-1}(\bm{y}^{k+1}-\bm{x}^{k+1}).
\end{equation*}
For compatible positive definite metrics, Assumption \ref{assmp:Bk} is equivalent in finite dimensions to the corresponding spectral and variation bounds in \cite{rockafellar2026advances}.

\item \textit{iPHA-BLP (relative-type).} For the fixed scalar metric $B_{k}=r^{-1}I$, $r>0$, we consider the method in \cite[Algorithm 1]{buscaglia2024inexact} with unit relaxation and the present multiplier convention. The scenario subsolvers provide candidate pairs with zero enlargement error. Collecting these pairs in $(\bm{y}^{k+1},\bm{u}^{k+1})$, we define
\begin{equation*}
\bm{u}^{k+1}\in\mathcal{T}(\bm{y}^{k+1}) \quad \text{and} \quad \bm{z}^{k+1}=\bm{x}^{k}+r^{-1}(\bm{w}^{k}-\bm{u}^{k+1}).
\end{equation*}
The candidate pair is accepted when
\begin{equation}\label{eq:blp-test}
\|\bm{y}^{k+1}-\bm{z}^{k+1}\|_{p}^{2} \leq \gamma^{2}\Bigl(\|\bm{x}^{k}-P_{\mathcal{N}}\bm{z}^{k+1} +P_{\mathcal{M}}\bm{y}^{k+1}\|_{p}^{2}+ \|\bm{x}^{k}-P_{\mathcal{N}}\bm{y}^{k+1} +P_{\mathcal{M}}\bm{z}^{k+1}\|_{p}^{2} \Bigr),
\end{equation}
where $0\leq\gamma<1$. The stopping tolerance depends on the current candidates $\bm{y}^{k+1}$ and $\bm{z}^{k+1}$, so it may change during the inner iterations.

Since $\bm{x}^{k}\in\mathcal{N}$ and $\mathcal{N}\perp \mathcal{M}$, the sum of the two squared norms on the right-hand side equals $\|\bm{z}^{k+1}-\bm{x}^{k}\|_{p}^{2} + \|\bm{y}^{k+1}-\bm{x}^{k}\|_{p}^{2}.$ Thus, the aggregate criterion is satisfied if each scenario candidate satisfies
\begin{equation*}
\|\bm{y}_{i}^{k+1}-\bm{z}_{i}^{k+1}\|^{2}\leq\gamma^{2}\big(\|\bm{z}_{i}^{k+1}-\bm{x}_{i}^{k}\|^{2}+\|\bm{y}_{i}^{k+1}-\bm{x}_{i}^{k}\|^{2}\big).
\end{equation*}
Once the candidate pair has been accepted, the primal-dual  update is
\begin{equation*}
\bm{x}^{k+1} = \bm{x}^{k}-\alpha_{k+1} \bigl(\bm{x}^{k}-P_{\mathcal{N}}\bm{z}^{k+1}\bigr), \qquad \bm{w}^{k+1} = \bm{w}^{k}-\alpha_{k+1}rP_{\mathcal{M}}\bm{y}^{k+1},
\end{equation*}
where
\begin{equation*}
\alpha_{k+1} = \frac{ \left\langle \bm{x}^{k}-P_{\mathcal{N}}\bm{z}^{k+1} +P_{\mathcal{M}}\bm{y}^{k+1}, \, \bm{x}^{k}-P_{\mathcal{N}}\bm{y}^{k+1} +P_{\mathcal{M}}\bm{z}^{k+1} \right\rangle_{p} }{ \|\bm{x}^{k}-P_{\mathcal{N}}\bm{z}^{k+1} +P_{\mathcal{M}}\bm{y}^{k+1}\|_{p}^{2} }.
\end{equation*}
If the denominator vanishes, the stopping criterion \eqref{eq:blp-test} with $\gamma<1$ gives $\bm{y}^{k+1}=\bm{z}^{k+1}=\bm{x}^{k}$. The definition of $\bm{z}^{k+1}$ then implies $\bm{u}^{k+1}=\bm{w}^{k}$. Since $\bm{u}^{k+1}\in\mathcal{T}(\bm{y}^{k+1})$, the current iterate is a solution. Unlike the variable-metric methods, iPHA-BLP uses a fixed scalar metric.

\item \textit{iParInv-BSS (relative-type).} The inexact method of partial inverse in \cite[Section 5, Equation (20)]{burachik2006inexact} can be applied to multistage SVIs. Each scenario subsolver provides a candidate pair $(\bm{y}_{i}^{k+1},\bm{u}_{i}^{k+1})$ and an enlargement error $\varepsilon_{i,k+1}\geq0$ satisfying $\bm{u}_{i}^{k+1}\in\mathcal{T}_{i}^{\varepsilon_{i,k+1}}(\bm{y}_{i}^{k+1})$. With unit scaling and unit proximal parameter, the candidate is accepted when
\begin{equation*}
\|\bm{y}_{i}^{k+1}+\bm{u}_{i}^{k+1}-\bm{x}_{i}^{k}-\bm{w}_{i}^{k}\|^{2} +2\varepsilon_{i,k+1} \leq \gamma^{2} \bigl( \|\bm{y}_{i}^{k+1}-\bm{x}_{i}^{k}\|^{2} + \|\bm{u}_{i}^{k+1}-\bm{w}_{i}^{k}\|^{2} \bigr),
\end{equation*}
where $0\leq\gamma<1/\sqrt{5}$ is the admissible range corresponding to $\widehat{\sigma}$ in \cite[Proposition 2.2 and Remarks 2.5, 3.1(2), and 4.1(2)]{costa2017inexact}. This criterion can be checked independently for each scenario, but its stopping tolerance depends on the current candidates and may change during the inner iterations. Once a candidate pair has been accepted for every scenario, the primal-dual update is
\begin{equation*}
\bm{x}^{k+1} = \bm{x}^{k}-P_{\mathcal{N}}\bm{u}^{k+1}, \qquad \bm{w}^{k+1} = \bm{w}^{k}-P_{\mathcal{M}}\bm{y}^{k+1}.
\end{equation*}
For $B_{k}=I$, this agrees with \eqref{eq:iPDA-update} upon substituting $\bm{\vartheta}^{k+1}=\bm{y}^{k+1}+\bm{u}^{k+1}-\bm{x}^{k}-\bm{w}^{k}$. For the fixed scalar metric $B_{k}=r^{-1}I$, we apply the method in coordinates obtained by scaling the primal and dual variables by $\sqrt{r}$ and $1/\sqrt{r}$, respectively; the numerical experiments use $r=10$.
\end{enumerate}

Like iParInv-BSS and iPHA-BLP, iVMPHA-Lag uses a fixed inexactness parameter without prescribing a summable tolerance sequence. While the stopping tolerances of these existing relative-type methods depend on current candidates, iVMPHA-Lag uses \eqref{eq:budget-update} to determine each scenario's tolerance before the inner iterations begin, keeping it fixed throughout these inner iterations.

%%%%%%%%%%%%%%%%%%%%%%%%%%%%%%%%%%%%%%%%%%%%%%%%%%%%%%%%%%%%
\section{Convergence analysis}\label{sec:convergence}

We establish global convergence of Algorithm \ref{algo:iPDA} through a descent estimate that combines the primal-dual distance with the error budget. Under metric subregularity at the limit point, we further prove local $R$-linear convergence of the primal-dual iterates. %The notation and metric conventions follow Section~\ref{sec:model}.

%%%%%%%%%%%%%%%%%%%%%%%%%%%%%%%%%%%%%%%%%%%%%%%%%%%%%%%%%%%%
\subsection{Global convergence}\label{subsec:global-convergence}

We begin by aggregating the scenario-wise quantities in Algorithm \ref{algo:iPDA}. Define
\begin{equation*}
\sigma_{k}^{2}=\sum_{i=1}^{s} p_{i}\sigma_{i,k}^{2},\quad s_{k+1}^{2}=\sum_{i=1}^{s} p_{i} s_{i,k+1}^{2},\quad \varepsilon_{k+1}=\sum_{i=1}^{s} p_{i}\varepsilon_{i,k+1},
\end{equation*}
and set
\begin{equation}\label{def:u}
\bm{u}^{k+1}=\bm{w}^{k}+\bm{\vartheta}^{k+1}-B_{k}^{-1}(\bm{y}^{k+1}-\bm{x}^{k}).
\end{equation}
The scenario conditions \eqref{eq:block-inexcond} and budget updates \eqref{eq:budget-update} imply the following aggregate relations:
\begin{subnumcases}{\label{eq:global-certificate}}
\;\bm{u}^{k+1}\in\mathcal{T}^{\varepsilon_{k+1}}(\bm{y}^{k+1}),\label{eq:global-certificate-a}\\
\;\|\bm{\vartheta}^{k+1}\|_{B_{k},p}^{2}+2\varepsilon_{k+1}\leq\rho^{2}\sigma_{k}^{2},\label{eq:global-certificate-b}\\
\;\sigma_{k+1}^{2}=s_{k+1}^{2}+\delta\sigma_{k}^{2}.\label{eq:global-certificate-c}
\end{subnumcases}
% \begin{subequations}\label{eq:global-certificate}
% \begin{empheq}[left=\empheqlbrace]{align}
% &\;\bm u^{k+1}\in\mathcal T^{\varepsilon_{k+1}}(\bm y^{k+1}),\label{eq:global-certificate-a}\\
% &\;\|\bm\vartheta^{k+1}\|_{B_k,p}^2+2\varepsilon_{k+1}\leq\rho^2\sigma_k^2,\label{eq:global-certificate-b}\\
% &\;\sigma_{k+1}^2=s_{k+1}^2+\delta\sigma_k^2.\label{eq:global-certificate-c}
% \end{empheq}
% \end{subequations}
Here, \eqref{eq:global-certificate-a} follows from the local-to-global enlargement relation \eqref{eq:weighted-enlargement}, while \eqref{eq:global-certificate-b}--\eqref{eq:global-certificate-c} follow by probability-weighted summation of the corresponding scenario relations. Combining \eqref{def:u} with the corrected update \eqref{eq:iPDA-update} and using Proposition \ref{prop:metric-identities}(iii) gives
\begin{equation}\label{eq:update-graph}
\bm{x}^{k+1}=\bm{x}^{k}-B_{k}P_{\mathcal{N}}\bm{u}^{k+1},\qquad \bm{w}^{k+1}=\bm{w}^{k}-B_{k}^{-1}P_{\mathcal{M}}\bm{y}^{k+1}.
\end{equation}
The aggregate progress measure in \eqref{eq:budget-update} can be written as
\begin{equation}\label{eq:s-graph}
s_{k+1}^{2}=\|P_{\mathcal{N}}(\bm{y}^{k+1}-\bm{x}^{k})\|_{B_{k}^{-1},p}^{2}+\|P_{\mathcal{M}}(\bm{u}^{k+1}-\bm{w}^{k})\|_{B_{k},p}^{2}.
\end{equation}
For later use, set $q:=\rho^{2}+\delta<1$. To obtain a descent estimate, we combine the squared primal-dual distance with the error budget. For any $(\bm{x}^{\star},\bm{w}^{\star})\in\Omega$, define
\begin{equation*}
D_{k}(\bm{x}^{\star},\bm{w}^{\star}):=\|(\bm{x}^{k}-\bm{x}^{\star},\bm{w}^{k}-\bm{w}^{\star})\|_{G_{k},p}^{2}\;,  \quad \Phi_{k}(\bm{x}^{\star},\bm{w}^{\star}):=D_{k}(\bm{x}^{\star},\bm{w}^{\star})+q\sigma_{k}^{2}.
\end{equation*}
The following lemma bounds the squared primal-dual distance and the updated error budget in terms of their current values. It then yields a quasi-Fej\'er estimate under the variable metric.

\begin{lemma}[One-step estimate]\label{lem:one-step}
Suppose that Assumption \ref{assmp:Bk} holds and $\Omega\neq\emptyset$. Let the sequence $\{({\bm{x}}^{k},{\bm{w}}^{k})\}$ be generated by Algorithm \ref{algo:iPDA}. Then, for every $k\geq0$ and every $(\bm{x}^{\star},\bm{w}^{\star})\in\Omega$,
\begin{equation}\label{eq:one-step}
\|(\bm{x}^{k+1}-\bm{x}^{\star},\bm{w}^{k+1}-\bm{w}^{\star})\|_{G_{k},p}^{2}+\sigma_{k+1}^{2}\leq D_{k}(\bm{x}^{\star},\bm{w}^{\star})+q\sigma_{k}^{2}.
\end{equation}
Consequently,
\begin{equation}\label{eq:quasi-fejer}
\Phi_{k+1}(\bm{x}^{\star},\bm{w}^{\star})\leq(1+\eta_{k})\Phi_{k}(\bm{x}^{\star},\bm{w}^{\star})-(1+\eta_{k}-q)\sigma_{k+1}^{2}.
\end{equation}
\end{lemma}

\begin{proof}
From \eqref{def:u} and Proposition \ref{prop:metric-identities}(iii), we have
\begin{equation}\label{eq:certificate-norm-expansion}
\|\bm{\vartheta}^{k+1}\|_{B_{k},p}^{2}=\|\bm{u}^{k+1}-\bm{w}^{k}\|_{B_{k},p}^{2}+\|\bm{y}^{k+1}-\bm{x}^{k}\|_{B_{k}^{-1},p}^{2}+2\langle\bm{y}^{k+1}-\bm{x}^{k},\bm{u}^{k+1}-\bm{w}^{k}\rangle_{p}.
\end{equation}
Subtracting \eqref{eq:s-graph} from \eqref{eq:certificate-norm-expansion} and using Proposition \ref{prop:metric-identities}(iii), together with $P_{\mathcal{N}}\bm{w}^{k}=0$ and $P_{\mathcal{M}}\bm{x}^{k}=0$, gives
\begin{equation}\label{eq:certificate-progress-balance}
\|\bm{\vartheta}^{k+1}\|_{B_{k},p}^{2}-s_{k+1}^{2}=\|P_{\mathcal{N}}\bm{u}^{k+1}\|_{B_{k},p}^{2}+\|P_{\mathcal{M}}\bm{y}^{k+1}\|_{B_{k}^{-1},p}^{2}+2\langle\bm{y}^{k+1}-\bm{x}^{k},\bm{u}^{k+1}-\bm{w}^{k}\rangle_{p}.
\end{equation}
Fix $(\bm{x}^{\star},\bm{w}^{\star})\in\Omega$. Then $\bm{w}^{\star}\in\mathcal{T}(\bm{x}^{\star})$, while \eqref{eq:global-certificate-a} gives $\bm{u}^{k+1}\in\mathcal{T}^{\varepsilon_{k+1}}(\bm{y}^{k+1})$. Hence, by the definition of the enlargement,
\begin{equation}\label{eq:enlargement-solution}
\langle\bm{y}^{k+1}-\bm{x}^{\star},\,\bm{u}^{k+1}-\bm{w}^{\star}\rangle_{p}\geq-\varepsilon_{k+1}.
\end{equation}

Next, using the update relations \eqref{eq:update-graph}, we obtain
\begin{equation}\label{eq:distance-increment}
\begin{aligned}
&\quad\;\|(\bm{x}^{k+1}-\bm{x}^{\star},\bm{w}^{k+1}-\bm{w}^{\star})\|_{G_{k},p}^{2}-D_{k}(\bm{x}^{\star},\bm{w}^{\star})\\
&= -2\langle\bm{x}^{k}-\bm{x}^{\star},\,P_{\mathcal{N}}\bm{u}^{k+1}\rangle_{p}-2\langle\bm{w}^{k}-\bm{w}^{\star},\,P_{\mathcal{M}}\bm{y}^{k+1}\rangle_{p}\\
&\quad\;+\|P_{\mathcal{N}}\bm{u}^{k+1}\|_{B_{k},p}^{2}+\|P_{\mathcal{M}}\bm{y}^{k+1}\|_{B_{k}^{-1},p}^{2} .
\end{aligned}
\end{equation}
Since $\bm{x}^{k}, \bm{x}^{\star}\in\mathcal{N}$ and $\bm{w}^{k}, \bm{w}^{\star}\in\mathcal{M}$, the two cross terms in \eqref{eq:distance-increment} satisfy
\begin{equation}\label{eq:cross-term-bound}
\begin{aligned}
&-\langle\bm{x}^{k}-\bm{x}^{\star},\,P_{\mathcal{N}}\bm{u}^{k+1}\rangle_{p}-\langle\bm{w}^{k}-\bm{w}^{\star},\,P_{\mathcal{M}}\bm{y}^{k+1}\rangle_{p}\\
={}&-\langle\bm{x}^{k}-\bm{x}^{\star},\,\bm{u}^{k+1}\rangle_{p}-\langle\bm{w}^{k}-\bm{w}^{\star},\,\bm{y}^{k+1}\rangle_{p}\\
={}&\;\langle\bm{y}^{k+1}-\bm{x}^{k},\,\bm{u}^{k+1}-\bm{w}^{k}\rangle_{p}-\langle\bm{y}^{k+1}-\bm{x}^{\star},\,\bm{u}^{k+1}-\bm{w}^{\star}\rangle_{p}\\
\leq{}&\;\langle\bm{y}^{k+1}-\bm{x}^{k},\,\bm{u}^{k+1}-\bm{w}^{k}\rangle_{p}+\varepsilon_{k+1},
\end{aligned}
\end{equation}
where the last inequality follows from \eqref{eq:enlargement-solution}. Substituting \eqref{eq:cross-term-bound} into \eqref{eq:distance-increment}, and then using \eqref{eq:certificate-progress-balance}, yields
\begin{equation}\label{eq:distance-certificate-bound}
\begin{aligned}
&\|(\bm{x}^{k+1}-\bm{x}^{\star}, \bm{w}^{k+1}-\bm{w}^{\star})\|_{G_{k},p}^{2}-D_{k}(\bm{x}^{\star},\bm{w}^{\star})\\
\leq{}&\;2\langle\bm{y}^{k+1}-\bm{x}^{k},\,\bm{u}^{k+1}-\bm{w}^{k}\rangle_{p}+2\varepsilon_{k+1}+\|P_{\mathcal{N}}\bm{u}^{k+1}\|_{B_{k},p}^{2}+\|P_{\mathcal{M}}\bm{y}^{k+1}\|_{B_{k}^{-1},p}^{2}\\
={}&\;\|\bm{\vartheta}^{k+1}\|_{B_{k},p}^{2}-s_{k+1}^{2}+2\varepsilon_{k+1} \leq \rho^{2}\sigma_{k}^{2}-s_{k+1}^{2},
\end{aligned}
\end{equation}
where the last inequality follows from \eqref{eq:global-certificate-b}. Using \eqref{eq:global-certificate-c} and $q=\rho^{2}+\delta$ in \eqref{eq:distance-certificate-bound} gives \eqref{eq:one-step}. Finally, by the successive-metric bound \eqref{eq:successive-metric},
\begin{align*}
\Phi_{k+1}(\bm{x}^{\star},\bm{w}^{\star})&\leq(1+\eta_{k})\|(\bm{x}^{k+1}-\bm{x}^{\star},\bm{w}^{k+1}-\bm{w}^{\star})\|_{G_{k},p}^{2}+q\sigma_{k+1}^{2}\\
&\leq(1+\eta_{k})\bigl(\Phi_{k}(\bm{x}^{\star},\bm{w}^{\star})-\sigma_{k+1}^{2}\bigr)+q\sigma_{k+1}^{2} \\
&=(1+\eta_{k})\Phi_{k}(\bm{x}^{\star},\bm{w}^{\star})-(1+\eta_{k}-q)\sigma_{k+1}^{2},
\end{align*}
where the second inequality follows from \eqref{eq:one-step}. This proves \eqref{eq:quasi-fejer}.
\end{proof}

%The same graph-variable representation also controls the actual primal-dual increment.

\begin{lemma}[Step bound]\label{lem:step}
Suppose that Assumption \ref{assmp:Bk} holds and $\Omega\neq\emptyset$. Let the sequence $\{({\bm{x}}^{k},{\bm{w}}^{k})\}$ be generated by Algorithm \ref{algo:iPDA}. Then, for every $k\geq0$,
\begin{equation}\label{eq:step-bound}
\|(\bm{x}^{k+1}-\bm{x}^{k},\bm{w}^{k+1}-\bm{w}^{k})\|_{G_{k},p}\leq s_{k+1}+\rho\sigma_{k} \leq \sigma_{k+1}+\rho\sigma_{k}.
\end{equation}
\end{lemma}

\begin{proof}
Using the corrected update \eqref{eq:iPDA-update}, the definition \eqref{def:u}, and the invariance of $\mathcal{N}$ and $\mathcal{M}$ under $B_{k}$ and $B_{k}^{-1}$, we obtain
\begin{equation*}
\bm{x}^{k+1}-\bm{x}^{k}=P_{\mathcal{N}}(\bm{y}^{k+1}-\bm{x}^{k})-B_{k}P_{\mathcal{N}}\bm{\vartheta}^{k+1}, \bm{w}^{k+1}-\bm{w}^{k}=P_{\mathcal{M}}(\bm{u}^{k+1}-\bm{w}^{k})-P_{\mathcal{M}}\bm{\vartheta}^{k+1}.
\end{equation*}
Hence,
\begin{equation*}
\begin{aligned}
&~~\|(\bm{x}^{k+1}-\bm{x}^{k},\bm{w}^{k+1}-\bm{w}^{k})\|_{G_{k},p}\\
\leq&~~ \|(P_{\mathcal{N}}(\bm{y}^{k+1}-\bm{x}^{k}), P_{\mathcal{M}}(\bm{u}^{k+1}-\bm{w}^{k}))\|_{G_{k},p} +\|(B_{k}P_{\mathcal{N}}\bm{\vartheta}^{k+1}, P_{\mathcal{M}}\bm{\vartheta}^{k+1})\|_{G_{k},p}\\
=&~~  s_{k+1}+\|\bm{\vartheta}^{k+1}\|_{B_{k},p} \leq s_{k+1}+\rho\sigma_{k},
\end{aligned}
\end{equation*}
where the equality follows from \eqref{eq:s-graph} and Proposition \ref{prop:metric-identities}(iii), and the last inequality follows from \eqref{eq:global-certificate-b}. Finally, \eqref{eq:global-certificate-c} gives $s_{k+1}\leq\sigma_{k+1}$, which proves \eqref{eq:step-bound}.
\end{proof}

This lemma shows that the changes between consecutive primal-dual iterates tend to zero whenever the error budgets tend to zero. Using the quasi-Fej\'er estimate and the inexactness condition, the following proposition establishes that the primal-dual iterates remain bounded, the budgets, residuals, and enlargement errors tend to zero, and the accepted candidates approach the current iterates. These properties will allow us to pass to the limit in the approximate inclusion and show that every cluster point is a solution.

\begin{proposition}\label{prop:asymptotic}
Suppose that Assumption \ref{assmp:Bk} holds and $\Omega\neq\emptyset$. Let the sequence $\{({\bm{x}}^{k},{\bm{w}}^{k})\}$ be generated by Algorithm \ref{algo:iPDA}. Then, for any $(\bm{x}^{\star},\bm{w}^{\star})\in\Omega$, the following statements hold.
\begin{itemize}
\item[(i)] The sequence $\{\Phi_{k}(\bm{x}^{\star},\bm{w}^{\star})\}$ converges to a finite limit, and the primal-dual sequence $\{(\bm{x}^{k},\bm{w}^{k})\}$ is bounded. Moreover,
\begin{equation}\label{eq:budget-sum}
\sum_{k=1}^{\infty}\sigma_{k}^{2}<\infty,\quad \sum_{k=0}^{\infty}\left(\|\bm{\vartheta}^{k+1}\|_{B_{k},p}^{2}+2\varepsilon_{k+1}\right)<\infty.
\end{equation}

\item[(ii)] The aggregate budgets, progress measures, residuals, and enlargement errors vanish:
\begin{equation*}
\sigma_{k}\to0,\quad s_{k}\to0,\quad \|\bm{\vartheta}^{k}\|_{p}\to0,\quad \varepsilon_{k}\to0.
\end{equation*}

\item[(iii)] The graph candidates are asymptotically consistent with the current primal-dual iterates:
\begin{equation}\label{eq:candidate-consistency}
\|\bm{y}^{k+1}-\bm{x}^{k}\|_{p}\to0,\quad \|\bm{u}^{k+1}-\bm{w}^{k}\|_{p}\to0.
\end{equation}

\item[(iv)] The primal-dual sequence is asymptotically regular:
\begin{equation}\label{eq:asymptotic-regularity}
\|(\bm{x}^{k+1}-\bm{x}^{k}, \bm{w}^{k+1}-\bm{w}^{k})\|_{p} \to0.
\end{equation}
\end{itemize}
\end{proposition}

\begin{proof}
\emph{Statement (i).} Fix $(\bm{x}^{\star},\bm{w}^{\star})\in\Omega$. By \eqref{eq:quasi-fejer} and the summability of $\{\eta_{k}\}$, \cite[Section~2.2, Lemma~2, p.~44]{p1987introduction} implies that $\{\Phi_{k}(\bm{x}^{\star},\bm{w}^{\star})\}$ converges to a finite limit. Hence there exists $\mu\geq0$ such that
\begin{equation*}
\lim_{k\to\infty} \Phi_{k}(\bm{x}^{\star},\bm{w}^{\star})=\mu.
\end{equation*}
In particular, $\{\Phi_{k}(\bm{x}^{\star},\bm{w}^{\star})\}$ is bounded. Since $D_{k}(\bm{x}^{\star},\bm{w}^{\star})\leq\Phi_{k}(\bm{x}^{\star},\bm{w}^{\star}),$ the uniform norm equivalence \eqref{eq:norm-equivalence} implies that $\{(\bm{x}^{k},\bm{w}^{k})\}$ is bounded.

Let $M>0$ be such that $\Phi_{k}(\bm{x}^{\star},\bm{w}^{\star})\leq M$ for all $k$. From \eqref{eq:quasi-fejer},
\begin{equation}\label{eq:budget-descent}
\begin{aligned}
(1-q)\sigma_{k+1}^{2}&\leq(1+\eta_{k}-q)\sigma_{k+1}^{2}\\
&\leq \Phi_{k}(\bm{x}^{\star},\bm{w}^{\star})-\Phi_{k+1}(\bm{x}^{\star},\bm{w}^{\star})+\eta_{k}\Phi_{k}(\bm{x}^{\star},\bm{w}^{\star})\\
&\leq \Phi_{k}(\bm{x}^{\star},\bm{w}^{\star})-\Phi_{k+1}(\bm{x}^{\star},\bm{w}^{\star})+M\eta_{k} .
\end{aligned}
\end{equation}
Summing \eqref{eq:budget-descent} from $k=0$ to $N$ gives
\begin{equation*}
(1-q)\sum_{k=0}^{N}\sigma_{k+1}^{2}\leq\Phi_{0}(\bm{x}^{\star},\bm{w}^{\star})-\Phi_{N+1}(\bm{x}^{\star},\bm{w}^{\star})+M\sum_{k=0}^{N}\eta_{k}.
\end{equation*}
Since $\Phi_{N+1}(\bm{x}^{\star},\bm{w}^{\star})\geq0$ and $\sum_{k=0}^{\infty}\eta_{k}<\infty$, letting $N\to\infty$ yields $\sum_{k=1}^{\infty}\sigma_{k}^{2}<\infty.$ Furthermore, by \eqref{eq:global-certificate-b},
\begin{equation*}
\sum_{k=0}^{\infty}\left(\|\bm{\vartheta}^{k+1}\|_{B_{k},p}^{2}+2\varepsilon_{k+1}\right)\leq\rho^{2}\sum_{k=0}^{\infty}\sigma_{k}^{2}<\infty,
\end{equation*}
where $\sigma_{0}<\infty$. This proves \eqref{eq:budget-sum}.

\medskip
\noindent
\emph{Statement (ii).} By \eqref{eq:budget-sum}, $\sigma_{k}\to0.$ From \eqref{eq:global-certificate-c}, $s_{k+1}^{2}\leq\sigma_{k+1}^{2},$ and hence $s_{k}\to0$. Moreover, \eqref{eq:global-certificate-b} gives
\begin{equation*}
0\leq\|\bm{\vartheta}^{k+1}\|_{B_{k},p}^{2}+2\varepsilon_{k+1}\leq\rho^{2}\sigma_{k}^{2}.
\end{equation*}
Together with $\sigma_{k}\to0$ and \eqref{eq:norm-equivalence}, this yields $\|\bm{\vartheta}^{k+1}\|_{p}\to0$ and $\varepsilon_{k+1}\to0$.

\medskip
\noindent
\emph{Statement (iii).} By \eqref{eq:s-graph} and $s_{k+1}\to0$, we have $\|P_{\mathcal{N}}(\bm{y}^{k+1}-\bm{x}^{k})\|_{B_{k}^{-1},p}\to0$ and $\|P_{\mathcal{M}}(\bm{u}^{k+1}-\bm{w}^{k})\|_{B_{k},p}\to0.$ The uniform norm equivalence therefore implies
\begin{equation}\label{eq:projected-consistency-p}
P_{\mathcal{N}}(\bm{y}^{k+1}-\bm{x}^{k})\to0, \quad P_{\mathcal{M}}(\bm{u}^{k+1}-\bm{w}^{k})\to0 \;\text{ in }\|\cdot\|_{p}.
\end{equation}

On the other hand, rearranging \eqref{def:u} gives
\begin{equation*}
\bm{\vartheta}^{k+1}=\bm{u}^{k+1}-\bm{w}^{k}+B_{k}^{-1}(\bm{y}^{k+1}-\bm{x}^{k}).
\end{equation*}
Projecting this identity onto $\mathcal{N}$ and $\mathcal{M}$, and using the invariance of these subspaces under $B_{k}^{-1}$, yields
\begin{align*}
P_{\mathcal{N}}(\bm{u}^{k+1}-\bm{w}^{k})&=P_{\mathcal{N}}\bm{\vartheta}^{k+1}-B_{k}^{-1}P_{\mathcal{N}}(\bm{y}^{k+1}-\bm{x}^{k}),\\
B_{k}^{-1}P_{\mathcal{M}}(\bm{y}^{k+1}-\bm{x}^{k})&=P_{\mathcal{M}}\bm{\vartheta}^{k+1}-P_{\mathcal{M}}(\bm{u}^{k+1}-\bm{w}^{k}).
\end{align*}
By Statement~(ii), \eqref{eq:projected-consistency-p}, and the uniform boundedness of $B_{k}$ and $B_{k}^{-1}$,
\begin{equation}\label{eq:complementary-consistency}
P_{\mathcal{N}}(\bm{u}^{k+1}-\bm{w}^{k})\to0,\quad P_{\mathcal{M}}(\bm{y}^{k+1}-\bm{x}^{k})\to0 \;\text{in }\|\cdot\|_{p}.
\end{equation}
Combining \eqref{eq:projected-consistency-p} and \eqref{eq:complementary-consistency} with $\mathcal{H}=\mathcal{N}\oplus\mathcal{M}$ proves \eqref{eq:candidate-consistency}.

\medskip
\noindent
\emph{Statement (iv).} By Lemma \ref{lem:step}, $\|(\bm{x}^{k+1}-\bm{x}^{k}, \bm{w}^{k+1}-\bm{w}^{k})\|_{G_{k},p}\leq s_{k+1}+\rho\sigma_{k}.$ Statement~(ii) shows that the right-hand side tends to zero. The uniform norm equivalence \eqref{eq:norm-equivalence} therefore gives \eqref{eq:asymptotic-regularity}.
\end{proof}

We are now ready to prove global convergence. We first show that every cluster point of the sequence $\{(\bm{x}^{k},\bm{w}^{k})\}$ generated by Algorithm \ref{algo:iPDA} solves the multistage SVI \eqref{SVI}. We then use the convergence of the potential $\Phi_{k}$ to show that the whole sequence converges to one such point.

\begin{theorem}[Global convergence]\label{thm:global}
Suppose that Assumption \ref{assmp:Bk} holds and $\Omega\neq\emptyset$. Then the primal-dual sequence $\{({\bm{x}}^{k},{\bm{w}}^{k})\}$ generated by Algorithm \ref{algo:iPDA} converges to some $(\bm{x}^{\infty},\bm{w}^{\infty})\in\Omega$. Moreover, $(\bm{y}^{k},\bm{u}^{k})$ also converges to $(\bm{x}^{\infty},\bm{w}^{\infty})$. In particular, $\bm{x}^{\infty}$ solves the multistage SVI \eqref{eq:svi-basic}.
\end{theorem}

\begin{proof}
By Proposition \ref{prop:asymptotic}(i), the sequence $\{(\bm{x}^{k},\bm{w}^{k})\}$ is bounded. Since $\mathcal{H}$ is finite-dimensional, it admits a cluster point. Let $(\bar{\bm{x}},\bar{\bm{w}})$ be one such cluster point, and let $\{k_{j}\}$ be a subsequence such that
\begin{equation}\label{eq:cluster-subsequence}
(\bm{x}^{k_{j}},\bm{w}^{k_{j}})\to(\bar{\bm{x}},\bar{\bm{w}}).
\end{equation}
Since every iterate belongs to $\mathcal{N}\times\mathcal{M}$ and both subspaces are closed, $(\bar{\bm{x}},\bar{\bm{w}})\in\mathcal{N}\times\mathcal{M}.$ By Proposition \ref{prop:asymptotic}(iii), we have $\bm{y}^{k+1}-\bm{x}^{k}\to0$ and $\bm{u}^{k+1}-\bm{w}^{k}\to0.$ Together with \eqref{eq:cluster-subsequence}, this gives
\begin{equation}\label{eq:candidate-cluster-limit}
\bm{y}^{k_{j}+1}\to\bar{\bm{x}}, \quad \bm{u}^{k_{j}+1}\to\bar{\bm{w}}.
\end{equation}

We next show that $\bar{\bm{w}}\in\mathcal{T}(\bar{\bm{x}})$. For any $(\bm{z},\bm{v})\in\operatorname{gph}\mathcal{T}$, \eqref{eq:global-certificate-a} and the definition of the enlargement give
\begin{equation}\label{eq:cluster-enlargement}
\langle \bm{y}^{k_{j}+1}-\bm{z},\, \bm{u}^{k_{j}+1}-\bm{v} \rangle_{p} \geq -\varepsilon_{k_{j}+1}.
\end{equation}
By Proposition \ref{prop:asymptotic}(ii), $\varepsilon_{k_{j}+1}\to0$. Passing to the limit in \eqref{eq:cluster-enlargement} using \eqref{eq:candidate-cluster-limit} yields
\begin{equation}\label{eq:cluster-monotone-relation}
\langle\bar{\bm{x}}-\bm{z},\,\bar{\bm{w}}-\bm{v}\rangle_{p}\geq0,\quad\forall\,(\bm{z},\bm{v})\in\operatorname{gph}\mathcal{T}.
\end{equation}
Since $\mathcal{T}$ is maximally monotone, \eqref{eq:cluster-monotone-relation} implies $\bar{\bm{w}}\in\mathcal{T}(\bar{\bm{x}})$.

It remains to show that the whole sequence converges to this cluster point. Since $(\bar{\bm{x}},\bar{\bm{w}})\in\Omega$, $(\bar{\bm{x}},\bar{\bm{w}})$ can be used as the reference solution in Proposition \ref{prop:asymptotic}(i). Hence, $\{\Phi_{k}(\bar{\bm{x}},\bar{\bm{w}})\}$ converges to a finite limit. On the other hand, \eqref{eq:cluster-subsequence} and the uniform norm equivalence imply $D_{k_{j}}(\bar{\bm{x}},\bar{\bm{w}})\to0.$ Since Proposition \ref{prop:asymptotic}(ii) also gives $\sigma_{k_{j}}\to0$, we obtain
\begin{equation*}
\Phi_{k_{j}}(\bar{\bm{x}},\bar{\bm{w}})=D_{k_{j}}(\bar{\bm{x}},\bar{\bm{w}})+q\sigma_{k_{j}}^{2}\to0.
\end{equation*}
Because the full sequence $\{\Phi_{k}(\bar{\bm{x}},\bar{\bm{w}})\}$ converges and one of its subsequences converges to zero, its limit must be zero. Thus, $\Phi_{k}(\bar{\bm{x}},\bar{\bm{w}})\to0.$ Combined with $0\leq D_{k}(\bar{\bm{x}},\bar{\bm{w}}) \leq \Phi_{k}(\bar{\bm{x}},\bar{\bm{w}}),$ and the uniform norm equivalence, we have $(\bm{x}^{k},\bm{w}^{k})\to(\bar{\bm{x}},\bar{\bm{w}}).$ Setting $(\bm{x}^{\infty},\bm{w}^{\infty}):=(\bar{\bm{x}},\bar{\bm{w}})$ proves the global convergence of $\{(\bm{x}^{k},\bm{w}^{k})\}$.

Finally, \eqref{eq:candidate-consistency} together with $(\bm{x}^{k},\bm{w}^{k})\to (\bm{x}^{\infty},\bm{w}^{\infty})$ gives $(\bm{y}^{k},\bm{u}^{k})\to (\bm{x}^{\infty},\bm{w}^{\infty})$. In the multistage SVI setting of Section \ref{sec:model}, $\bm{x}^{\infty}$ therefore solves \eqref{eq:svi-basic}.
\end{proof}

%%%%%%%%%%%%%%%%%%%%%%%%%%%%%%%%%%%%%%%%%%%%%%%%%%%%%%%%%%%%
\subsection{Linear convergence}\label{subsec:linear-convergence}

To establish a local linear rate, we use the following primal-dual residual mapping. For $(\bm{x},\bm{w})\in\mathcal{N}\times\mathcal{M}$, define
\begin{equation}\label{eq:residual-map}
\mathcal{R}_{\mathcal{N}}(\bm{x},\bm{w}) := \{\bm{a}+\bm{c}: \bm{a}\in\mathcal{N},\ \bm{c}\in\mathcal{M},\ \bm{a}+\bm{w}\in\mathcal{T}(\bm{x}+\bm{c})\}.
\end{equation}
By the definition of the partial inverse \eqref{eq:partial-inverse}, we have $\mathcal{R}_{\mathcal{N}}(\bm{x},\bm{w})=\mathcal{T}_{\mathcal{N}}(\bm{x}+\bm{w}).$ Moreover, the isometry and zero-set correspondence in Section \ref{subsec:partial-inverse} give
\begin{equation}\label{eq:EB-isometry}
\begin{aligned}
\operatorname{dist}_{p}((\bm{x},\bm{w}),\Omega) &= \operatorname{dist}_{p}( \bm{x}+\bm{w},\mathcal{T}_{\mathcal{N}}^{-1}(0)),  \\
\operatorname{dist}_{p}(0,\mathcal{R}_{\mathcal{N}}(\bm{x},\bm{w})) &= \operatorname{dist}_{p}( 0,\mathcal{T}_{\mathcal{N}}(\bm{x}+\bm{w})).
\end{aligned}
\end{equation}

We impose the following local primal-dual error bound at the limit point obtained in Theorem \ref{thm:global}.

\begin{assumption}[Local primal-dual error bound]\label{assmp:EB}
Let $(\bm{x}^{\infty},\bm{w}^{\infty})$ be the limit in Theorem \ref{thm:global}. There exist $r_{0}>0$ and $\kappa>0$ such that
\begin{equation*}
\operatorname{dist}_{p}((\bm{x},\bm{w}),\Omega) \leq \kappa\operatorname{dist}_{p}( 0,\mathcal{R}_{\mathcal{N}}(\bm{x},\bm{w}))
\end{equation*}
for every $(\bm{x},\bm{w})\in\mathcal{N}\times\mathcal{M}$ satisfying $\|(\bm{x}-\bm{x}^{\infty},\bm{w}-\bm{w}^{\infty})\|_{p}<r_{0}.$
\end{assumption}

\begin{remark}[Relation to metric subregularity]\label{rem:EB-equivalence}
Let $\bm{z}^{\infty}:=\bm{x}^{\infty}+\bm{w}^{\infty}.$ By $\mathcal{R}_{\mathcal{N}}(\bm{x},\bm{w})=\mathcal{T}_{\mathcal{N}}(\bm{x}+\bm{w})$ and \eqref{eq:EB-isometry}, Assumption \ref{assmp:EB} is equivalent to the metric subregularity of $\mathcal{T}_{\mathcal{N}}$ at $(\bm{z}^{\infty},\bm{0})$:
\begin{equation*}
\operatorname{dist}_{p}( \bm{z},\mathcal{T}_{\mathcal{N}}^{-1}(0)) \leq \kappa\operatorname{dist}_{p}( 0,\mathcal{T}_{\mathcal{N}}(\bm{z}))
\end{equation*}
for all $\bm{z}$ sufficiently close to $\bm{z}^{\infty}$. The isometry preserves both the constant $\kappa$ and the neighborhood radius. In particular, no isolated solution or uniform error bound around the whole solution set is required. For affine $F_{i}$ and polyhedral $C_{i}$, the operator $\mathcal{T}$ is piecewise polyhedral, and so is its partial inverse $\mathcal{T}_{\mathcal{N}}$: its graph is a finite union of polyhedra, a property preserved by the linear graph transformation defining the partial inverse. Assumption \ref{assmp:EB} follows from the local error-bound property of such multifunctions \cite{robinson1981continuity}.
\end{remark}

The next estimate connects this local error bound to the progress and error budget. An exact regularized solution is used only as a comparison point in the proof.

\begin{proposition}[Distance estimate]\label{prop:distance-estimate}
Suppose that Assumption \ref{assmp:Bk} and Assumption \ref{assmp:EB} hold and $\Omega\neq\emptyset$. Let the sequence $\{({\bm{x}}^{k},{\bm{w}}^{k})\}$ be generated by Algorithm \ref{algo:iPDA}. Set $\tau:=\Lambda(1+\kappa).$ Then, for all sufficiently large $k$,
\begin{equation}\label{eq:distance-estimate}
\operatorname{dist}_{G_{k+1},p} \bigl((\bm{x}^{k+1},\bm{w}^{k+1}),\Omega\bigr) \leq \tau(s_{k+1}+\rho\sigma_{k}) \leq \tau(\sigma_{k+1}+\rho\sigma_{k}).
\end{equation}
\end{proposition}

\begin{proof}
Consider the exact regularized inclusion
\begin{equation*}
0\in \mathcal{T}(\bm{y})-\bm{w}^{k}+B_{k}^{-1}(\bm{y}-\bm{x}^{k}).
\end{equation*}
Since $\mathcal{T}+B_{k}^{-1}(\cdot-\bm{x}^{k})-\bm{w}^{k}$ is maximally and strongly monotone, it has a unique zero. Let $(\widetilde{\bm{y}},\widetilde{\bm{u}})$ denote the corresponding exact graph pair:
\begin{equation}\label{eq:exact-comparison-pair}
\widetilde{\bm{u}}\in\mathcal{T}(\widetilde{\bm{y}}), \qquad \widetilde{\bm{u}} = \bm{w}^{k}-B_{k}^{-1}(\widetilde{\bm{y}}-\bm{x}^{k}).
\end{equation}

From \eqref{def:u} and \eqref{eq:exact-comparison-pair}, we have $\bm{\vartheta}^{k+1} = \bm{u}^{k+1}-\widetilde{\bm{u}} + B_{k}^{-1}(\bm{y}^{k+1}-\widetilde{\bm{y}}).$ Since $\bm{u}^{k+1}\in \mathcal{T}^{\varepsilon_{k+1}}(\bm{y}^{k+1})$ and $\widetilde{\bm{u}}\in\mathcal{T}(\widetilde{\bm{y}})$, the enlargement inequality gives
\begin{equation*}
\langle \bm{y}^{k+1}-\widetilde{\bm{y}},\, \bm{u}^{k+1}-\widetilde{\bm{u}} \rangle_{p} \geq -\varepsilon_{k+1}.
\end{equation*}
Expanding the squared norm of $\bm{\vartheta}^{k+1}$ therefore yields
\begin{equation}\label{eq:resolvent-comparison}
\|\bm{y}^{k+1}-\widetilde{\bm{y}}\|_{B_{k}^{-1},p}^{2} + \|\bm{u}^{k+1}-\widetilde{\bm{u}}\|_{B_{k},p}^{2}\leq \|\bm{\vartheta}^{k+1}\|_{B_{k},p}^{2} +2\varepsilon_{k+1} \leq \rho^{2}\sigma_{k}^{2},
\end{equation}
where the last inequality follows from \eqref{eq:global-certificate-b}.

Define $\widetilde{\bm{x}}:=P_{\mathcal{N}}\widetilde{\bm{y}}$ and $\widetilde{\bm{w}}:=P_{\mathcal{M}}\widetilde{\bm{u}}.$ Taking $\bm{a}=P_{\mathcal{N}}\widetilde{\bm{u}}$ and $\bm{c}=P_{\mathcal{M}}\widetilde{\bm{y}}$ in \eqref{eq:residual-map} gives
\begin{equation}\label{eq:comparison-residual}
P_{\mathcal{N}}\widetilde{\bm{u}} + P_{\mathcal{M}}\widetilde{\bm{y}} \in \mathcal{R}_{\mathcal{N}}(\widetilde{\bm{x}},\widetilde{\bm{w}}).
\end{equation}
Moreover, by \eqref{eq:s-graph}, \eqref{eq:resolvent-comparison}, and the metric orthogonal decompositions, we obtain
\begin{equation}\label{eq:comparison-step}
\|(\widetilde{\bm{x}}-\bm{x}^{k}, \widetilde{\bm{w}}-\bm{w}^{k})\|_{G_{k},p} \leq s_{k+1}+\rho\sigma_{k}.
\end{equation}

Projecting \eqref{eq:exact-comparison-pair} onto $\mathcal{N}$ and $\mathcal{M}$ gives $P_{\mathcal{N}}\widetilde{\bm{u}}+P_{\mathcal{M}}\widetilde{\bm{y}}= B_{k}^{-1}(\bm{x}^{k}-\widetilde{\bm{x}})+B_{k}(\bm{w}^{k}-\widetilde{\bm{w}}).$ Hence, by \eqref{eq:comparison-step} and the uniform norm equivalence, we have
\begin{equation}\label{eq:comparison-residual-bound}
\|P_{\mathcal{N}}\widetilde{\bm{u}} + P_{\mathcal{M}}\widetilde{\bm{y}}\|_{p} \leq \sqrt{\Lambda}\, (s_{k+1}+\rho\sigma_{k}).
\end{equation}

By Theorem \ref{thm:global}, Proposition \ref{prop:asymptotic}(ii), and \eqref{eq:comparison-step}, we have $(\widetilde{\bm{x}},\widetilde{\bm{w}}) \to (\bm{x}^{\infty},\bm{w}^{\infty}).$ Therefore, Assumption \ref{assmp:EB} applies to $(\widetilde{\bm{x}},\widetilde{\bm{w}})$ for all sufficiently large $k$. Using \eqref{eq:comparison-residual} and \eqref{eq:comparison-residual-bound}, we obtain
\begin{equation}\label{eq:comparison-error-bound}
\operatorname{dist}_{p} \bigl((\widetilde{\bm{x}},\widetilde{\bm{w}}),\Omega\bigr) \leq \kappa\sqrt{\Lambda}\, (s_{k+1}+\rho\sigma_{k}).
\end{equation}

Moreover, from \eqref{eq:update-graph} and \eqref{eq:exact-comparison-pair}, we have
\begin{equation}\label{eq:iterate-comparison-bound}
\begin{aligned}
&\quad~\|(\bm{x}^{k+1}-\widetilde{\bm{x}}, \bm{w}^{k+1}-\widetilde{\bm{w}})\|_{G_{k},p}^{2}  \\
&= \|P_{\mathcal{N}}(\bm{u}^{k+1}-\widetilde{\bm{u}})\|_{B_{k},p}^{2} + \|P_{\mathcal{M}}(\bm{y}^{k+1}-\widetilde{\bm{y}})\|_{B_{k}^{-1},p}^{2}\\
&\leq \|\bm{u}^{k+1}-\widetilde{\bm{u}}\|_{B_{k},p}^{2} + \|\bm{y}^{k+1}-\widetilde{\bm{y}}\|_{B_{k}^{-1},p}^{2} \leq \rho^{2}\sigma_{k}^{2},
\end{aligned}
\end{equation}
where the last inequality follows from \eqref{eq:resolvent-comparison}.

Finally, using the base product norm as an intermediate norm, \eqref{eq:comparison-error-bound}, \eqref{eq:iterate-comparison-bound}, and \eqref{eq:distance-equivalence}, we obtain, for all sufficiently large $k$,
\begin{equation*}
\begin{aligned}
\operatorname{dist}_{G_{k+1},p} \bigl((\bm{x}^{k+1},\bm{w}^{k+1}),\Omega\bigr) &\leq \sqrt{\Lambda}\, \operatorname{dist}_{p} \bigl((\bm{x}^{k+1},\bm{w}^{k+1}),\Omega\bigr)\\
&\leq \sqrt{\Lambda}\left[ \operatorname{dist}_{p} \bigl((\widetilde{\bm{x}},\widetilde{\bm{w}}),\Omega\bigr) + \|(\bm{x}^{k+1}-\widetilde{\bm{x}}, \bm{w}^{k+1}-\widetilde{\bm{w}})\|_{p} \right]\\
&\leq \Lambda\kappa(s_{k+1}+\rho\sigma_{k}) + \Lambda\rho\sigma_{k}\\
&\leq \Lambda(1+\kappa)(s_{k+1}+\rho\sigma_{k}).
\end{aligned}
\end{equation*}
Thus the first inequality in \eqref{eq:distance-estimate} holds with $\tau=\Lambda(1+\kappa)$. The second follows from $s_{k+1}\leq\sigma_{k+1}$ by \eqref{eq:global-certificate-c}.
\end{proof}

The distance estimate and the one-step descent inequality now yield a contraction of an augmented potential. The resulting decay of the budgets, together with the step bound in Lemma \ref{lem:step}, then gives a linear rate for the actual primal-dual iterates.

\begin{theorem}[Local linear convergence]\label{thm:linear}
Suppose that Assumption \ref{assmp:Bk} and Assumption \ref{assmp:EB} hold and $\Omega\neq\emptyset$. Let the sequence $\{({\bm{x}}^{k},{\bm{w}}^{k})\}$ be generated by Algorithm \ref{algo:iPDA}. Choose $\nu\in(0,1-q)$ and define
\begin{equation}\label{eq:linear-merit}
\Psi_{k}:=\operatorname{dist}_{G_{k},p}^{2}\bigl((\bm{x}^{k},\bm{w}^{k}),\Omega\bigr)+(q+\nu)\sigma_{k}^{2}.
\end{equation}
Let $\ell:=\min\{\nu,1-q-\nu\}$, $L:=\tau^{2}(1+\rho)+q+\nu$, and $\mu:=\frac{L}{L+\ell}\in(0,1)$. Then,
\begin{equation}\label{eq:linear-contraction}
\Psi_{k+1} \leq \mu(1+\eta_{k})\Psi_{k},
\end{equation}
for all sufficiently large $k$. Moreover, the following statements hold.
\begin{enumerate}[label=(\roman*),leftmargin=*]
\item $\{\Psi_{k}\}$ converges to zero at a $Q$-linear rate.
%For every $\widehat\mu\in(\mu,1)$, there exists $K_{\widehat\mu}$
%such that
%\(
%\Psi_{k+1}\leq\widehat\mu\,\Psi_k,
%~
%k\geq K_{\widehat\mu}.
%\)
%Hence $\{\Psi_k\}$ converges to zero $Q$-linearly.

\item $\{(\bm{x}^{k},\bm{w}^{k})\}$ converges to $(\bm{x}^{\infty},\bm{w}^{\infty})$ $R$-linearly with $\|(\bm{x}^{k}-\bm{x}^{\infty}, \bm{w}^{k}-\bm{w}^{\infty})\|_{p}= \mathcal{O}(\mu^{k/2}).$
%\begin{equation}\label{eq:iterate-rate}
%\|(\bm x^k-\bm x^\infty,
%\bm w^k-\bm w^\infty)\|_p=
%\mathcal O(\mu^{k/2}).
%\end{equation}
\end{enumerate}
\end{theorem}

\begin{proof}
Taking the infimum over $(\bm{x}^{\star},\bm{w}^{\star})\in\Omega$ in \eqref{eq:one-step} and then applying the successive-metric bound \eqref{eq:successive-metric}, we obtain
\begin{equation}\label{eq:Psi-descent}
\begin{aligned}
\Psi_{k+1} &\leq (1+\eta_{k})\Psi_{k} -(1+\eta_{k}-q-\nu)\sigma_{k+1}^{2} -(1+\eta_{k})\nu\sigma_{k}^{2}\\
&\leq (1+\eta_{k})\Psi_{k} -\ell(\sigma_{k+1}^{2}+\sigma_{k}^{2}),
\end{aligned}
\end{equation}
where $\ell>0$ because $\nu\in(0,1-q)$.

On the other hand, Proposition \ref{prop:distance-estimate} gives, for all sufficiently large $k$, that
\begin{equation}\label{eq:Psi-budget-bound}
\begin{aligned}
\Psi_{k+1} &\leq \tau^{2}(\sigma_{k+1}+\rho\sigma_{k})^{2} +(q+\nu)\sigma_{k+1}^{2}  \\
&\leq \tau^{2}(1+\rho) \bigl(\sigma_{k+1}^{2}+\rho\sigma_{k}^{2}\bigr) +(q+\nu)\sigma_{k+1}^{2} \leq L(\sigma_{k+1}^{2}+\sigma_{k}^{2}),
\end{aligned}
\end{equation}
where we used $(a+\rho b)^{2}\leq(1+\rho)(a^{2}+\rho b^{2})$ and $0\leq\rho<1$. Substituting \eqref{eq:Psi-budget-bound} into \eqref{eq:Psi-descent} yields
\begin{equation*}
\Psi_{k+1}\leq(1+\eta_{k})\Psi_{k}-\frac{\ell}{L}\Psi_{k+1},
\end{equation*}
and hence
\begin{equation*}
\Psi_{k+1}\leq\frac{L}{L+\ell}(1+\eta_{k})\Psi_{k}=\mu(1+\eta_{k})\Psi_{k}.
\end{equation*}

Since $\eta_{k}\to0$, for every $\widehat{\mu}\in(\mu,1)$ there exists $K_{\widehat{\mu}}$ such that $\mu(1+\eta_{k})\leq\widehat{\mu}, ~ k\geq K_{\widehat{\mu}}.$ Statement~(i) follows.

To prove Statement~(ii), fix $K$ sufficiently large that \eqref{eq:linear-contraction} holds for all $k\geq K$. Iterating \eqref{eq:linear-contraction} gives
\begin{equation}\label{eq:Psi-geometric-bound}
\Psi_{k} \leq \Psi_{K}\mu^{k-K} \prod_{j=K}^{k-1}(1+\eta_{j}) \leq C\mu^{k-K}, \qquad k\geq K,
\end{equation}
for some $C>0$, since $\prod_{j=K}^{\infty}(1+\eta_{j})\leq\exp\!\left(\sum_{j=K}^{\infty}\eta_{j}\right)<\infty.$ Because $q+\nu>0$, \eqref{eq:linear-merit} and \eqref{eq:Psi-geometric-bound} imply $\sigma_{k}=\mathcal{O}(\mu^{k/2}).$ Moreover, \eqref{eq:global-certificate-c} gives $s_{k+1}\leq\sigma_{k+1}$. Therefore Lemma \ref{lem:step} and \eqref{eq:norm-equivalence} yield
\begin{equation}\label{eq:increment-linear-rate}
\|(\bm{x}^{k+1}-\bm{x}^{k}, \bm{w}^{k+1}-\bm{w}^{k})\|_{p} = \mathcal{O}(\mu^{k/2}).
\end{equation}
Using Theorem \ref{thm:global} and summing the tail of \eqref{eq:increment-linear-rate}, we obtain
\begin{equation*}
\|(\bm{x}^{k}-\bm{x}^{\infty}, \bm{w}^{k}-\bm{w}^{\infty})\|_{p} \leq \sum_{j=k}^{\infty} \|(\bm{x}^{j+1}-\bm{x}^{j}, \bm{w}^{j+1}-\bm{w}^{j})\|_{p} = \mathcal{O}(\mu^{k/2}).
\end{equation*}
\end{proof}

\begin{remark}[Rate parameters and the exact case]\label{rem:rate-parameters}
A balanced choice is $\nu=\frac{1-q}{2},$ for which
\begin{equation*}
\ell=\frac{1-q}{2}\quad \text{and} \quad \mu = \frac{2\tau^{2}(1+\rho)+1+q}{2\tau^{2}(1+\rho)+2} <1.
\end{equation*}
This factor is a proof-based upper bound for the contraction of the augmented potential and is not intended to be sharp. In particular, when $\rho=0$, the algorithm solves each scenario subproblem exactly and the iterates are independent of the retained budget, whereas the bound above still depends on $q=\delta$ through the augmented potential.

For variable metrics, the one-step factor in \eqref{eq:linear-contraction} is $\mu(1+\eta_{k})$, so any $\widehat{\mu}\in(\mu,1)$ is an eventual fixed $Q$-linear factor. When $\eta_{k}\equiv0$, the bound reduces directly to $\mu$. The primal-dual iterates retain the $R$-linear order $\mathcal{O}(\mu^{k/2})$ because the cumulative metric variation is finite.
\end{remark}

%%%%%%%%%%%%%%%%%%%%%%%%%%%%%%%%%%%%%%%%%%%%%%%%%%%%%%%%%%%%
\section{Numerical experiments}\label{sec:numerical}

All computations are performed in MATLAB R2024b on Windows 11 using a PC equipped with an Intel Core i7-12700 processor and 16 GB of memory. All runs are executed serially without parallel workers.

%%%%%%%%%%%%%%%%%%%%%%%%%%%%%%%%%%%%%%%%%%%%%%%%%%%%%%%%%%%%
\subsection{Test problem generation}\label{subsec:test-problems}

We generate synthetic multistage SVIs on complete $m$-ary scenario trees with $T$ stages and $s=m^{T-1}$ leaf scenarios. For every scenario $i$, the feasible set is the box $C_{i}=[0,1]^{d}$, where $d=Td_{0}$ and $d_{0}$ is the dimension of the decision vector at each stage. At each nonterminal node, independent weights drawn uniformly from $[0.5,1.5]$ are assigned to its outgoing branches and normalized to obtain their conditional probabilities. The probability $p_{i}$ of scenario $i$ is the product of the conditional branch probabilities along its root-to-leaf path.

The scenario operators are generated from the common form
\begin{equation*}
F_{i}(\bm{x}_{i})=(H_{i}+S_{i})\bm{x}_{i}+\zeta D_{i}\left(\bm{x}_{i}-\tfrac{1}{2}\bm{1}\right)^{\circ 3}+\bm{q}_{i},
\end{equation*}
where $\bm{1}\in\mathbb{R}^{d}$ is the all-ones vector and $(\cdot)^{\circ 3}$ denotes the componentwise cubic power. Let $\omega_{t}=100^{(t-1)/(T-1)}$ for $t=1,\ldots,T$. The diagonal entries of $D_{i}$ corresponding to stage $t$ are $(D_{i})_{\ell\ell}=\omega_{t}\exp(\xi_{i\ell})$, where the $\xi_{i\ell}$ are drawn independently and uniformly from $[-0.6,0.6]$. Thus, the nominal stage-scale contrast is $\omega_{T}/\omega_{1}=100$. The symmetric and skew-symmetric matrices are defined by
\begin{equation*}
H_{i}=D_{i}^{1/2}\left(0.05I+L_{i}^{\top} L_{i}\right)D_{i}^{1/2}, \quad S_{i}=\kappa D_{i}^{1/2}K_{i}D_{i}^{1/2},
\end{equation*}
where $L_{i}$ is a cyclic difference matrix with independent edge weights drawn uniformly from $[0.5,1.5]$, normalized so that $\lambda_{\max}(L_{i}^{\top} L_{i})=1$. The matrix $K_{i}$ is obtained by skew-symmetrizing a Gaussian random matrix and normalizing it to satisfy $\|K_{i}\|_{2}=1$. Since $H_{i}\succ0$, $S_{i}^{\top}=-S_{i}$, and the cubic term is monotone for $\zeta\geq0$, each $F_{i}$ is strongly monotone.

To construct a known solution $(\bm{x}^{\star},\bm{w}^{\star})\in\Omega$, we draw the decision vectors at each node uniformly from $\{0,0.25,0.5,0.75,1\}$, assigning the same values to scenarios sharing that node, and obtain $\bm{w}^{\star}$ by projecting a Gaussian vector onto $\mathcal{M}$. For each scenario $i$, we choose $\bm{v}_{i}^{\star}\in N_{C_{i}}(\bm{x}_{i}^{\star})$ and determine $\bm{q}_{i}$ from $\bm{w}_{i}^{\star}=F_{i}(\bm{x}_{i}^{\star})+\bm{v}_{i}^{\star}$. The constructed solution is used for instance generation and validation but is not provided to the algorithms.

We consider four problem regimes specified by $(\zeta,\kappa)\in\{0,1\}\times\{0.25,1\}$, where $\zeta=0$ and $\zeta=1$ correspond to affine and nonlinear problems, respectively. To examine the effect of problem size, we vary one of $T$, $m$, and $d_{0}$ at a time from the baseline $(T,m,d_{0})=(4,4,20)$, as summarized in Table \ref{tab:test-sizes}. This gives seven distinct structural configurations. For each configuration and regime, we generate ten independent instances using prescribed random seeds, with the same instances used by all methods.

\begin{table}[!htbp]
\centering
\caption{Problem-size sweeps used in the numerical experiments.}\label{tab:test-sizes}
\begin{tabular}{lcccc}
\hline
Sweep & Varied parameter & Fixed parameters & $s$ & $d$\\
\hline
Stage
& $T\in\{3,4,5\}$
& $(m,d_{0})=(4,20)$
& $16,64,256$
& $60,80,100$\\
Branch
& $m\in\{2,4,6\}$
& $(T,d_{0})=(4,20)$
& $8,64,216$
& $80$\\
Dimension
& $d_{0}\in\{10,20,50\}$
& $(T,m)=(4,4)$
& $64$
& $40,80,200$\\
\hline
\end{tabular}
\end{table}

%%%%%%%%%%%%%%%%%%%%%%%%%%%%%%%%%%%%%%%%%%%%%%%%%%%%%%%%%%%%
\subsection{Comparisons with existing methods}\label{subsec:method-comparison}

We compare iVMPHA-Lag with the exact PHA, the relative-type methods iPHA-BLP \cite{buscaglia2024inexact} and iParInv-BSS \cite{burachik2006inexact}, and Rockafellar's absolute-type methods iVMPDA-absI and iVMPDA-absII \cite{rockafellar2026advances}. We use the scalar metric $B_{i,k}=0.1I$ for exact PHA, iPHA-BLP, iParInv-BSS, and iVMPHA-Lag (denoted as iVMPHA-Lag-S), and the same adaptive metric update rule for iVMPDA-absI, iVMPDA-absII, and iVMPHA-Lag (denoted as iVMPHA-Lag-A).

All methods start from $(\bm{x}^{0},\bm{w}^{0})=(\tfrac{1}{2}\bm{1},\bm{0})$ and use the common stopping criterion
\begin{equation*}
\operatorname{res}(\bm{x},\bm{w}):= \frac{\|\bm{x}-P_{\mathcal{C}}(\bm{x}-\mathcal{F}(\bm{x})+\bm{w})\|_{p}}{1+\|\bm{x}\|_{p}}\leq 10^{-6}.
\end{equation*}
Here, $\mathcal{C}=\prod_{i=1}^{s}C_{i}$ and $\mathcal{F}(\bm{x})=(F_{i}(\bm{x}_{i}))_{i=1}^{s}$, as in Section \ref{subsec:primal-dual}, and $P_{\mathcal{C}}$ denotes the projection onto $\mathcal{C}$. The residual is evaluated initially and after each outer iteration. Runs are limited to $5000$ outer iterations and $5000$ inner iterations per scenario subproblem; those exceeding either limit without meeting the corresponding stopping criterion are classified as unsuccessful.

% The scalar configurations use $R_{k,i}:=B_{k,i}^{-1}=10I$ for all $i$ and $k$. For the adaptive configurations, we initialize $R_{0,i}=10I$ and use the stage-wise regularization matrix $R_{k,i}=\operatorname{Diag}(r_{1,k}I_{d_0},\ldots,r_{T,k}I_{d_0})$, with coefficients shared by all scenarios. At outer iterations $K\in\{5,10,20,40\}$, we estimate the stage-wise operator scales using the secant ratios
% $$c_{t,K}=\left(\frac{\sum_{j\in\mathcal J_K}\sum_i p_i\|\Delta F_{i,t}^j\|^2}
% {\sum_{j\in\mathcal J_K}\sum_i p_i\|\Delta\bm y_{i,t}^j\|^2}\right)^{1/2},$$
% where $\Delta F_{i,t}^j$ and $\Delta\bm y_{i,t}^j$ denote the changes in the stage-$t$ components of the operator values and accepted candidates between consecutive outer iterations, and $\mathcal J_K=\{\max(2,K-4),\ldots,K\}$. The estimated scales determine provisional regularization coefficients
% \[
% \widetilde r_{t,K}=\frac{10c_{t,K}}{\left(\prod_{\ell=1}^{T}c_{\ell,K}\right)^{1/T}},\quad t=1,\ldots,T.
% \]
% Thus, the coefficients are proportional to the estimated operator scales, with geometric mean $10$. Their logarithmic deviations from $10$ are uniformly reduced, if necessary, to obtain target coefficients $r_{t,K}\in[1,100]$. The current coefficients are then updated toward these targets in logarithmic scale, preserving their geometric mean and limiting each change to a factor of two. The metric remains fixed after $K=40$. Since the stage coefficients are common to all scenarios, uniformly bounded, and updated only finitely many times, the resulting metric sequence satisfies Assumption~\ref{assmp:Bk}.

For iVMPDA-absI, iVMPDA-absII, and iVMPHA-Lag-A, we construct a variable-metric sequence $\{B_{k}\}$ from the scalar initialization $B_{0}=0.1I$ using the following stage-wise adaptive rule. Specifically, we use the same metric in every scenario $i=1,\ldots,s$, and set $B_{i,k}$ as \eqref{stage-metric}, so that $B_{k}\mathcal{N}\subseteq \mathcal{N}$. Let $K$ denote the number of completed outer iterations. Metric updates are permitted only when $K\in\{5,10,20,40\}$ and the metric remains fixed otherwise. At each update, we first estimate the stage-wise operator scales using the secant ratios
\begin{equation*}
\chi_{t,K}:=\left(\frac{\sum_{j\in\mathcal{J}_{K}}\sum_{i} p_{i}\|\Delta F_{i,t}^{j}\|^{2}}{\sum_{j\in\mathcal{J}_{K}}\sum_{i} p_{i}\|\Delta\bm{y}_{i,t}^{j}\|^{2}}\right)^{1/2}, \quad \text{with} \quad \mathcal{J}_{K}:=\{\max(2,K-4),\ldots,K\},
\end{equation*}
where $\Delta F_{i,t}^{j}$ and $\Delta\bm{y}_{i,t}^{j}$ denote successive differences in the stage-$t$ components of $F_{i}(\bm{y}_{i}^{j})$ and $\bm{y}_{i}^{j}$, respectively. Let $h_{t,K-1}=-\log(10b_{t,K-1})$ and $\bm{h}_{K-1}=(h_{t,K-1})_{t=1}^{T}$. Let $g_{t,K}=\log\chi_{t,K}-T^{-1}\sum_{\ell=1}^{T}\log\chi_{\ell,K}$ and $\bm{g}_{K}=(g_{t,K})_{t=1}^{T}$. The metric coefficients are then updated by
\begin{equation*}
\bm{h}_{K}^{*}=\frac{\bm{g}_{K}}{\max\{1,\|\bm{g}_{K}\|_{\infty}/\log 10\}}, \; \bm{h}_{K}=\bm{h}_{K-1}+\alpha_{K}(\bm{h}_{K}^{*}-\bm{h}_{K-1}), \; b_{t,K}=0.1e^{-h_{t,K}},
\end{equation*}
where $\alpha_{K}=\min\left\{1,\frac{\log 2}{\|\bm{h}_{K}^{*}-\bm{h}_{K-1}\|_{\infty}}\right\},$ with $\alpha_{K}=1$ when the denominator vanishes. Thus, larger secant ratios $\chi_{t,K}$ correspond to smaller target metric coefficients, toward which the current coefficients are updated while ensuring $(\prod_{t=1}^{T} b_{t,K})^{1/T}=0.1$, $0.01\leq b_{t,K}\leq 1$, and $1/2\leq b_{t,K}/b_{t,K-1}\leq 2$.
Updates with nearly zero secant denominators or invalid scale estimates are skipped, and the metric remains unchanged at all other iterations. Consequently, the resulting metric sequence satisfies Assumption \ref{assmp:Bk}.

For iVMPDA-absI and iVMPDA-absII, we prescribe the summable tolerance sequences $a_{k}=10^{-3}(k+1)^{-1.1}$ and $c_{k}=0.5(k+1)^{-1.1}$, respectively. We set $\gamma_{\mathrm{BSS}}=0.4$ for iParInv-BSS, $\gamma_{\mathrm{BLP}}=0.6$ for iPHA-BLP, and $(\rho,\delta)=(0.6,0.16)$ for iVMPHA-Lag, whose initial error budgets are $\sigma_{i,0}=\max\{10^{-3},\|F_{i}(\bm{x}_{i}^{0})-\bm{w}_{i}^{0}\|_{B_{i,0}}\}.$ The exact PHA requires $\|\bm{\vartheta}_{i}^{k+1}\|_{B_{i,k}}\leq 10^{-10}$ for each scenario. All parameters are fixed across problem instances.

All scenario subproblems are solved using the same diagonally preconditioned extragradient method \cite{korpelevich1976extragradient}, with preconditioner $Q_{i,k}$ and stepsize $\beta_{i,k}$ given by
\begin{equation*}
Q_{i,k} =\operatorname{Diag}(\operatorname{diag}H_{i}) +0.75\,\zeta D_{i}+B_{i,k}^{-1},\quad  \beta_{i,k} =0.8/L_{i,k},
\end{equation*}
where $L_{i,k}=(\|\widehat{J}_{i,k}\|_{1}\|\widehat{J}_{i,k}\|_{\infty})^{1/2}$ and $\widehat{J}_{i,k}=Q_{i,k}^{-1/2}(H_{i}+S_{i}+B_{i,k}^{-1}+0.75\,\zeta D_{i})Q_{i,k}^{-1/2}$. The inner iterations are initialized at the previously accepted candidates. We set $\varepsilon_{i,k+1}=0$ and evaluate the residuals according to \eqref{eq:normal-certificate}. Each inexact method follows its acceptance criterion in Section \ref{sec:iPDA}: iVMPHA-Lag and iParInv-BSS use scenario-wise tests, whereas iVMPDA-absI, iVMPDA-absII, and iPHA-BLP use aggregate tests.

Table \ref{tab:baseline-comparison} compares the six methods under seven configurations on the baseline problem $(T,m,d_{0})=(4,4,20)$. The reported outer iterations, scenario evaluations $E$, and elapsed times are averaged over ten instances for each problem regime. All runs are successful. The count $E$ includes all evaluations of the scenario operators $F_{i}$, including those used for initialization and outer stopping tests. It is shown that all inexact methods substantially reduce the computational cost relative to exact PHA. Under the common scalar metric, iVMPHA-Lag-S performs comparably to iParInv-BSS and is more efficient than iPHA-BLP. Under the common adaptive metric update rule, iVMPHA-Lag-A requires fewer scenario evaluations and less time than iVMPDA-absI and iVMPDA-absII in all four regimes. In particular, its outer-iteration counts are comparable to those of iVMPDA-absII, while it requires only about one quarter of the scenario evaluations, indicating substantially less inner computational work.

\begin{table}[!htbp]
\centering
\caption{Comparison on the baseline configuration $(T,m,d_{0})=(4,4,20)$}\label{tab:baseline-comparison}
\small
\setlength{\tabcolsep}{3pt}
\renewcommand{\arraystretch}{1.0}
\begin{tabular}{@{}clcrr@{\hspace{6pt}}crr@{}}
\toprule
&
& \multicolumn{3}{c}{Affine ($\zeta=0$)}
& \multicolumn{3}{c}{Nonlinear ($\zeta=1$)} \\
\cmidrule(r{5pt}){3-5}
\cmidrule(l{5pt}){6-8}
& Method
& OutIter & $E\quad\;$ & Time(s)
& OutIter & $E\quad\;$ & Time(s) \\
\midrule

\multicolumn{1}{c|}{%
\multirow{7}{*}{\shortstack[c]{Lower\\Coupling\\$(\kappa=0.25)$}}}
& PHA
& 610  & 8,464,744 & 165.02
& 368  & 5,725,125 & 137.05 \\

\multicolumn{1}{c|}{}
& iPHA-BLP
& 1319 & 592,070  & 11.62
& 889  & 289,082  & 7.29 \\

\multicolumn{1}{c|}{}
& iParInv-BSS
& 607  & 508,176  & 10.16
& 368  & 260,910  & 5.98 \\

\multicolumn{1}{c|}{}
& iVMPHA-Lag-S
& 607 & 519,625 & 10.28
& 368 & 254,490 & 5.68 \\

\multicolumn{1}{c|}{}
& iVMPDA-absI
& 446  & 368,544  & 6.98
& 478  & 274,643  & 6.51 \\

\multicolumn{1}{c|}{}
& iVMPDA-absII
& 133  & 419,846  & 7.71
& 86   & 273,216  & 5.72 \\

\multicolumn{1}{c|}{}
& iVMPHA-Lag-A
& 129
& 102,906
& 2.14
& 82
& 63,204
& 1.62 \\

\midrule

\multicolumn{1}{c|}{%
\multirow{7}{*}{\shortstack[c]{Higher\\Coupling\\$(\kappa=1)$}}}
& PHA
& 345 & 8,662,506 & 171.74
& 268 & 6,202,082 & 141.55 \\

\multicolumn{1}{c|}{}
& iPHA-BLP
& 814 & 681,709 & 13.85
& 576 & 344,179 & 8.01 \\

\multicolumn{1}{c|}{}
& iParInv-BSS
& 347 & 622,980 & 13.51
& 280 & 323,583 & 8.42 \\

\multicolumn{1}{c|}{}
& iVMPHA-Lag-S
& 347 & 622,856 & 13.06
& 282 & 310,570 & 7.45 \\

\multicolumn{1}{c|}{}
& iVMPDA-absI
& 422 & 519,386 & 9.69
& 469 & 405,043 & 8.76 \\

\multicolumn{1}{c|}{}
& iVMPDA-absII
& 83  & 500,525 & 8.93
& 66  & 354,298 & 7.42 \\

\multicolumn{1}{c|}{}
& iVMPHA-Lag-A
& 88
& 122,720
& 2.67
& 69
& 79,573
& 1.95 \\

\bottomrule
\end{tabular}
\end{table}

\begin{figure}[!htbp]
\centering
\includegraphics[width=0.9\textwidth]{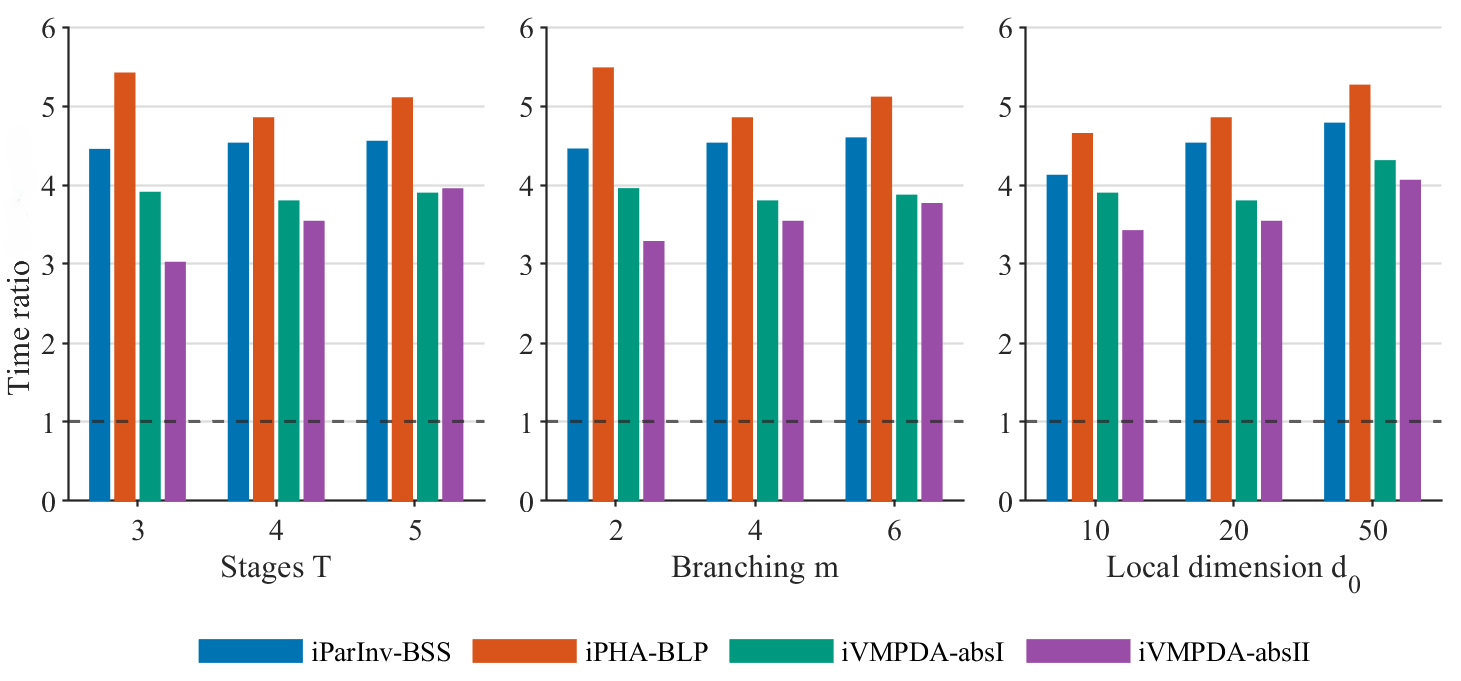}
\caption{Relative elapsed times across the three problem-size sweeps in Table \ref{tab:test-sizes}. Each bar shows the ratio of the competing method's mean elapsed time to that of iVMPHA-Lag-A over the same 40 instances, with 10 instances from each problem regime.}\label{fig:scal-comp}
\end{figure}

Figure \ref{fig:scal-comp} compares the elapsed times of the four competing inexact methods relative to iVMPHA-Lag-A across the three problem-size sweeps in $T$, $m$, and $d_{0}$. At each setting of $(T,m,d_0)$, each bar gives the ratio of the competing method's mean elapsed time to that of iVMPHA-Lag-A on the same $40$ instances ($10$ per regime). These ratios range from about $3$ to $5.5$, which shows that iVMPHA-Lag-A maintains an average time advantage over these methods across the tested problem sizes.

%%%%%%%%%%%%%%%%%%%%%%%%%%%%%%%%%%%%%%%%%%%%%%%%%%%%%%%%%%%%
\subsection{Influence of the error-budget design}\label{subsec:budget-influence}

We examine the effects of lagged error control and positive retention by comparing three variants of iVMPHA-Lag under the common fixed stage-scaled metric $B_{i,k}:=\operatorname{Diag}(\omega_{1}^{-1}I_{d_{0}},\ldots,\omega_{T}^{-1}I_{d_{0}})$ for all $i,k$. The diagnostic variant Current-0 replaces \eqref{eq:block-inexcond-b} in Step 1 of Algorithm \ref{algo:iPDA} with
\begin{equation*}
\|\bm{\vartheta}_{i}^{k+1}\|_{B_{i,k}}^{2}+2\varepsilon_{i,k+1}\leq\rho^{2}\sigma_{i,k+1}^{2},
\end{equation*}
where $\sigma_{i,k+1}^{2}=s_{i,k+1}^{2}$ is determined by the current candidates with $\delta=0$. Since this budget depends on coordinated information from all scenarios, its evaluation requires a global check. We implement Current-0 synchronously, with all scenarios continuing refinement after a rejected check. In contrast, Lagged-0 and Lagged-0.16 use the preassigned budgets of Algorithm \ref{algo:iPDA}, with $\delta=0$ and $\delta=0.16$, respectively, allowing independent scenario-wise stopping. All three variants use $\rho=0.6$, the same primal-dual initial point, inner solver, and outer stopping criterion. We test them on ten common instances in each of the four problem regimes with $(T,m,d_{0})=(5,4,20)$.

Table \ref{tab:budget-influence} reports the results, where all runs are successful. Here, ``RejAtmps'' counts the failed inner stopping checks that require further refinement, while ``ActSn/Round'' measures the average number of scenarios continuing refinement after each rejected check. The three variants require nearly identical numbers of outer iterations. Compared with Current-0, Lagged-0 substantially reduces the number of active scenarios per refinement round, together with scenario evaluations and elapsed time. Positive retention further reduces rejected attempts and scenario evaluations, with relatively little change in the number of active scenarios per round. These results indicate that lagged stopping and positive retention reduce inner computational work without materially affecting outer progress.

\begin{table}[!htbp]
\centering
\caption{Influence of error-budget design on iVMPHA-Lag for $(T,m,d_{0})=(5,4,20)$.}\label{tab:budget-influence}
\small
\setlength{\tabcolsep}{2.5pt}
\renewcommand{\arraystretch}{0.95}
\begin{tabular}{l!{\hspace{8pt}}lcrr!{\hspace{6pt}}cc}
\toprule
($\zeta,\kappa$) & Variant & OutIter & $E\quad\;$ & Time(s) & RejAtmps & \shortstack{ActSn/Round} \\
\midrule
\multirow{3}{*}{\shortstack[l]{($0,0.25$)}}
& Current-0   & 195 & 748,595 & 31.51 & 1364 & 256 \\
& Lagged-0    & 194 & 646,683 & 28.01 & 1357 & 220 \\
& Lagged-0.16 & 194 & 625,131 & 27.89 & 1339 & 215 \\
\midrule
\multirow{3}{*}{\shortstack[l]{($0,1$)}}
& Current-0   & 131 & 915,021 & 34.19 & 1721 & 256 \\
& Lagged-0    & 130 & 737,625 & 29.23 & 1636 & 215 \\
& Lagged-0.16 & 129 & 694,204 & 27.38 & 1533 & 215 \\
\midrule
\multirow{3}{*}{\shortstack[l]{($1,0.25$)}}
& Current-0   & 134 & 585,062 & 26.84 & 1075 & 256 \\
& Lagged-0    & 133 & 481,244 & 23.35 & 1058 & 211 \\
& Lagged-0.16 & 133 & 456,635 & 22.25 & 1021 & 207 \\
\midrule
\multirow{3}{*}{\shortstack[l]{($1,1$)}}
& Current-0   & 103 & 711,398 & 27.76 & 1337 & 256 \\
& Lagged-0    & 103 & 560,847 & 22.76 & 1251 & 213 \\
& Lagged-0.16 & 102 & 518,554 & 20.67 & 1171 & 210 \\
\bottomrule
\end{tabular}
\end{table}

%%%%%%%%%%%%%%%%%%%%%%%%%%%%%%%%%%%%%%%%%%%%%%%%%%%%%%%%%%%%
\section{Conclusion}\label{sec:conclusion}

We propose iVMPHA-Lag, a variable-metric inexact PHA for finite-scenario multistage monotone SVIs. Its lagged scenario budgets allow independent inner stopping with tolerances fixed before the inner iterations begin, without prescribing a summable tolerance sequence. Under the stated assumptions, we establish finite inner termination and global primal-dual convergence, with local $R$-linear convergence under metric subregularity at the limit point. Numerical experiments on synthetic strongly monotone box-constrained SVIs demonstrate the computational benefits of adaptive metrics, lagged stopping, and positive budget retention. Future work may address total inner work for specific subsolvers and parallel implementations.

\appendix
% Unique appendix hyperlink destinations for the installed LaTeX engine.
\renewcommand{\theHsection}{\Alph{section}}
\renewcommand{\theHtheorem}{\theHsection.\arabic{theorem}}
\renewcommand{\theHequation}{\theHsection.\arabic{equation}}

%%%%%%%%%%%%%%%%%%%%%%%%%%%%%%%%%%%%%%%%%%%%%%%%%%%%%%%%%%%%
\section{Correspondence with a VM-HPE error condition}\label{app:vm-hpe}

We relate the probability-weighted aggregate iteration to the VM-HPE framework \cite{goncalves2018variable}. The primal-dual iterate is represented by $\bm{z}^{k}$, the accepted scenario pair by $\widehat{\bm{z}}^{k+1}$, and the aggregate error budget by an auxiliary nonnegative scalar $e_{k}$:
\begin{equation*}
\bm{z}^{k}=\bm{x}^{k}+\bm{w}^{k},\; \widehat{\bm{z}}^{k+1}=P_{\mathcal{N}}\bm{y}^{k+1}+P_{\mathcal{M}}\bm{u}^{k+1},\; M_{k}=B_{k}^{-1}P_{\mathcal{N}}+B_{k}P_{\mathcal{M}},\; e_{k}=q\sigma_{k}^{2}.
\end{equation*}
Here $q=\rho^{2}+\delta$, as in Section \ref{subsec:global-convergence}. The scalar $e_{k}$ is a multiple of the aggregate error budget, distinct from the enlargement error $\varepsilon_{k+1}$. The operators $M_{k}$ are positive definite and self-adjoint in $\mathcal{H}$. Their restrictions to ${{\mathcal{N}}}$ and ${\mathcal{M}}$ are $B_{k}^{-1}$ and $B_{k}$, respectively, so Assumption \ref{assmp:Bk} implies $\frac{1}{1+\eta_{k}}M_{k}\preceq M_{k+1}\preceq(1+\eta_{k})M_{k}.$ Equations \eqref{eq:update-graph} and \eqref{eq:s-graph} give
\begin{equation}\label{eq:hpe-identities}
\begin{aligned}
&\quad\quad\;\bm{r}^{k+1}:=M_{k}(\bm{z}^{k}-\bm{z}^{k+1}) =P_{\mathcal{N}}\bm{u}^{k+1}+P_{\mathcal{M}}\bm{y}^{k+1},\\
&\|\bm{z}^{k+1}-\widehat{\bm{z}}^{k+1}\|_{M_{k},p}^{2} =\|\bm{\vartheta}^{k+1}\|_{B_{k},p}^{2},\quad \|\bm{z}^{k}-\widehat{\bm{z}}^{k+1}\|_{M_{k},p}^{2}=s_{k+1}^{2}.
\end{aligned}
\end{equation}
By the partial-inverse enlargement identity \eqref{eq:partial-enlargement}, \eqref{eq:global-certificate} implies
\begin{equation}\label{eq:hpe-extended}
\bm{r}^{k+1}\in \mathcal{T}_{\mathcal{N}}^{\varepsilon_{k+1}}(\widehat{\bm{z}}^{k+1}),~ \|\bm{z}^{k+1}-\widehat{\bm{z}}^{k+1}\|_{M_{k},p}^{2}+2\varepsilon_{k+1}+e_{k+1} \leq q\|\bm{z}^{k}-\widehat{\bm{z}}^{k+1}\|_{M_{k},p}^{2}+e_{k}.
\end{equation}
Indeed, the left side of the second inequality is at most $\rho^{2}\sigma_{k}^{2}+q(s_{k+1}^{2}+\delta \sigma_{k}^{2})$. Since $q=\rho^{2}+\delta<1$, one has $q-(\rho^{2}+q\delta)=\delta(1-q)\geq0$. The inclusion uses \eqref{eq:partial-enlargement} at the candidate in \eqref{eq:hpe-identities}.

When $\varepsilon_{k+1}=0$ at every iteration, \eqref{eq:hpe-extended} is the VM-HPE inclusion and error condition in \cite{goncalves2018variable}, with its parameter $\sigma$ being set to $q$ and its auxiliary error set to $e_{k}$, after shifting the metric index. This framework parameter $\sigma$ is distinct from the square roots $\sigma_{i,k}$ of our scenario error budgets. Algorithm \ref{algo:iPDA} meets this aggregate condition through scenario stopping criteria whose tolerances are fixed before the inner iterations begin. Section \ref{sec:convergence} also allows nonzero enlargement errors and proves local $R$-linear convergence of the primal-dual iterates under Assumption \ref{assmp:EB}.

%%%%%%%%%%%%%%%%%%%%%%%%%%%%%%%%%%%%%%%%%%%%%%%%%%%%%%%%%%%%
% References use a local copy of the full manuscript's BibTeX database.
\bibliographystyle{plain}
\bibliography{Ref_iVMPHA_Lag}

\end{document}